\documentclass[11pt]{article}

\usepackage{hyperref,amsmath,amsfonts,latexsym,url,color,subfig,enumerate,graphics,graphicx}
\usepackage{longtable}
\usepackage{mathabx}
\usepackage{mathtools}
\usepackage{wasysym}

\newcommand{\m}[1]{{\bf{#1}}}
\newcommand{\g}[1]{\boldsymbol #1}

\newcommand{\C}[1]{{\cal {#1}}}

\newcommand{\lsup}[2]{{\vphantom{#2}}^{#1}{{#2}}}
\newcommand{\lrsup}[3]{{\vphantom{#2}}^{#1}{{#2}}^{#3}}
\newcommand{\lsupddt}[2]{{\vphantom{\frac{d#2}{dt}}}^{#1}{\frac{d#2}{dt}}}

  \title{{\bf Minimum-Time Earth-to-Mars Interplanetary Orbit Transfer Using Adaptive Gaussian Quadrature Collocation}}
  \author{Brittanny V.~Holden\footnote{Ph.D. Student, Department of Mechanical and Aerospace Engineering, University of Florida.  E-mail:  brittannyholden@ufl.edu} \\ Shan He\footnote{Undergraduate Student, Department of Mechanical and Aerospace Engineering, University of Florida.  E-mail:  shanhe0824@ufl.edu} \\ Anil V.~Rao\footnote{Professor, Department of Mechanical and Aerospace Engineering, University of Florida.  Associate Fellow, AIAA.  E-mail: anilvrao@ufl.edu.  Corresponding Author.} \vspace{12pt} \\ {\em University of Florida} \\ {\em Gainesville, FL 32611}}  
\date{}

\begin{document}

\maketitle
  
\begin{abstract}
The problem of minimum-time, low-thrust, Earth-to-Mars interplanetary orbital trajectory optimization is considered.  The minimum-time orbital transfer problem is modeled as a four-phase optimal control problem where the four phases correspond to planetary alignment, Earth escape, heliocentric transfer, and Mars capture.  The four-phase optimal control problem is then solved using a direct collocation adaptive Gaussian quadrature collocation method.  The following three models are used in the study: (1) circular planetary motion; (2) elliptic planetary motion; and (3) elliptic planetary motion with gravity perturbations, where the transfer begins in a geostationary orbit and terminates in a Mars-stationary orbit. Results for all three cases are provided, and one particular case is studied in detail to show the key features of the optimal solutions. Using the particular value  thrust specific force of $9.8\times 10^{-4}~\textrm{m}\cdot\textrm{s}^{-2}$, it was found that the minimum times for cases (1), (2), and (3) are, respectively, 215 d, 196 d, and 198 d with departure dates, respectively, of 1 July 2020, 30 June 2020, and 28 June 2020.  Finally, the problem formulation developed in this study is compared against prior work on an Earth-to-Mars interplanetary orbit transfer where it is found that the results of this research show significant improvement in transfer time relative to the prior work.
\end{abstract}

\section*{Nomenclature\label{sect:nomenclature}}

\renewcommand{\baselinestretch}{1}
\normalsize\normalfont
\begin{longtable}{ll}

$a$ & thrust specific force  \vspace{2pt}\\
$a_{pr}$ & radial component of gravity perturbations  \vspace{2pt}\\
$a_{p\theta}$ & transverse component of gravity perturbations  \vspace{2pt}\\
$D$ & length unit \vspace{2pt}\\
$e_E$ & eccentricity of Earth \vspace{2pt}\\
$\textbf{e}_E$ & eccentricity vector of Earth \vspace{2pt}\\
$e_M$ & eccentricity of Mars \vspace{2pt}\\
$\textbf{e}_M$ & eccentricity vector of Mars \vspace{2pt}\\
$f_E$ & second modified equinoctial orbital element of Earth  \vspace{2pt}\\
$f_M$ & second modified equinoctial orbital element of Mars \vspace{2pt}\\ 
$g_E$ & third modified equinoctial orbital element of Earth \vspace{2pt}\\
$g_M$ & third modified equinoctial orbital element of Mars \vspace{2pt}\\
$L_E$ & true longitude of Earth \vspace{2pt}\\
$L_M$ & true longitude of Mars \vspace{2pt}\\
$L_{M_0}$ & initial true longitude value of Mars \vspace{2pt}\\ 
$p_E$ & semi-latus rectum of Earth \vspace{2pt}\\
$p_M$ & semi-latus rectum of Mars \vspace{2pt}\\ 
$P$ & phase number \vspace{2pt}\\ 
$r$ & distance from central body to spacecraft \vspace{2pt}\\
$r_E$ & distance from Earth to Sun \vspace{2pt}\\
$r_M$ & distance from Mars to Sun \vspace{2pt}\\
$r_{pe}$ & distance from the spacecraft to Earth in the heliocentric transfer phase \vspace{2pt}\\
$r_{pm}$ & distance from the spacecraft to Mars in the heliocentric transfer phase \vspace{2pt}\\
$r_{psE}$ & distance from the spacecraft to Sun in the Earth escape phase \vspace{2pt}\\
$r_{psM}$ & distance from the spacecraft to Sun in the Mars capture phase \vspace{2pt}\\
$R_E$ & radius of Earth \vspace{2pt}\\
$R_M$ & radius of Mars \vspace{2pt}\\
$R_E^{SOI}$ & sphere of influence of Earth $\vspace{2pt}$\\
$R_M^{SOI}$ & sphere of influence of Mars $\vspace{2pt}$\\
$R_{SE}$ & distance from Sun to Earth \vspace{2pt}\\
$R_{SM}$ & distance from Sun to Mars \vspace{2pt}\\
$t$ & time \vspace{2pt}\\
$t_0$ & initial time \vspace{2pt}\\ 
$t_f$ & terminal time \vspace{2pt}\\ 
$T$ & time unit \vspace{2pt}\\
$v_r$ & radial component of spacecraft velocity \vspace{2pt}\\
$v_\theta$ & transverse component of spacecraft velocity \vspace{2pt}\\
$V$ & speed unit \vspace{2pt}\\
$w_r$ & radial component of thrust direction \vspace{2pt}\\
$w_\theta$ & transverse component of thrust direction \vspace{2pt}\\  
$\theta$  & longitude of spacecraft measured from the line of Aries \vspace{2pt}\\  
$\mu$ & gravitational parameter of central body \vspace{2pt}\\
$\mu_E$ & gravitational parameter of Earth \vspace{2pt}\\
$\mu_S$ & gravitational parameter of Sun \vspace{2pt}\\
$\mu_M$ & gravitational parameter of Mars \vspace{2pt}\\
$\nu_E$ & true anomaly of Earth \vspace{2pt}\\
$\nu_M$ & true anomaly of Mars \vspace{2pt}\\ 
$\varpi_E$ & longitude of perihelion of Earth orbit \vspace{2pt}\\
$\varpi_M$ & longitude of perihelion of Mars orbit \vspace{2pt}\\
$\Aries$ & line of Aries \vspace{2pt}\\
$\cancer$ & line of Cancer \vspace{2pt}\\
\end{longtable}

\renewcommand{\baselinestretch}{2}\normalsize\normalfont

\section{Introduction}\label{Introduction}

Interplanetary space travel has been a topic of interest to the space community for several decades.  Due to the large distances between planets, most interplanetary orbital transfers require the expenditure of a large amount of propellant.  Traditionally, large spacecraft that used high-thrust chemical propulsion were used for interplanetary travel.  More recently, high-thrust chemical propulsion has been replaced with low-thrust propulsion using either electric propulsion or solar electric propulsion. The main benefit of using a low-thrust propulsion system is that the fuel consumption is significantly lower when compared with the fuel consumption using high-thrust chemical propulsion.  One drawback to using low-thrust propulsion over high-thrust propulsion is that an orbital transfer using low-thrust propulsion is significantly longer than an orbital transfer using high-thrust propulsion due to the fact that a low-thrust engine produces significantly less thrust than a high-thrust engine.  This research focuses on the use of low-thrust propulsion to accomplish an interplanetary orbital transfer from Earth to Mars.  

A great deal of research has been done previously on the design of interplanetary orbital transfers \cite{Rayman2002,ByrnesLonguskiAldrin1993,Dargent2004,Breakwell1966,NahVadaliBraden2001,Petukhov2012,EnglanderConway2017,Yam2010, WallConway2005, GondelachNoomen2015,WuWangPohXu2009,D'AmarioByrnesStanford1982, Miele1997, TangConway1995}.    Refs.~\cite{Rayman2002, ByrnesLonguskiAldrin1993} both study specific problems to further expand the knowledge of low-thrust spacecraft and interplanetary orbital transfers. Ref.~\cite{Rayman2002} studied the Deep Space 1 mission because it was the first interplanetary mission to be propelled by a low-thrust solar electric propulsion system.  A key result of Ref.~\cite{Rayman2002} is that constraints on the spacecraft attitude and periods of coasting or thrusting of the propulsion system cause various trajectory design issues.  Ref.~\cite{ByrnesLonguskiAldrin1993} studied the Aldrin orbit, which is the simplest cycler orbit between Earth and Mars, and determines the minimum-impulse optimal solution from a multiconic optimization program for a 15 year cycle.  Refs.~\cite{Dargent2004,Breakwell1966, NahVadaliBraden2001,Petukhov2012} solved minimum-fuel optimal control problems by solving the Hamiltonian boundary-value problem (HBVP) arising from the calculus of variations.  In particular, Ref.~\cite{NahVadaliBraden2001} used an indirect multiple-shooting method to solve a three-dimensional minimum-fuel Earth-to-Mars orbital trajectory for a low-thrust spacecraft, while Ref.~\cite{Petukhov2012} transformed the HBVP into the Cauchy problem through the use of the continuation (homotopic) method.  In a manner similar to that used in Ref.~\cite{Petukhov2012},  Refs.~\cite{EnglanderConway2017,Yam2010,WallConway2005} employed a hybrid optimization method in order to solve their desired interplanetary orbital trajectory optimal control problem by implementing genetic algorithms, basin hopping, calculus of variations, and/or an adaptive neighborhood global optimization algorithm. Specifically, Ref.~\cite{WallConway2005} solved a three-phase minimum-fuel Earth-to-Mars orbital trajectory optimal control problem by creating a novel algorithm that forced the spacecraft to arrive at Mars and then optimized the final mass of the fuel. Refs.~\cite{GondelachNoomen2015, WuWangPohXu2009, D'AmarioByrnesStanford1982} developed novel methods to solve orbital transfer optimal control problems involving low-thrust spacecraft. Ref.~\cite{GondelachNoomen2015} shaped the velocity components as a function of time and polar angle during the transfer, Ref.~\cite{Herman2002} showed the robustness of a higher order collocation 7th degree system, Ref.~\cite{WuWangPohXu2009} utilized the Legendre pseudospectral method and SNOPT, a sparse nonlinear optimization software, to solve a minimum-fuel satellite formation maneuver, and Ref.~\cite{D'AmarioByrnesStanford1982} solved constrained interplanetary trajectory optimization problems by developing a computer program, PLATO (Planetary Trajectory Optimization). Ref.~\cite{Miele1997} solved a four-phase interplanetary orbital transfer of a spacecraft from a low Earth orbit to a low Martian orbit by utilizing a variable-stepsize integration technique implemented with the sequential gradient-restoration algorithm. Ref.~\cite{TangConway1995} numerically solved a minimum-time, low-thrust interplanetary transfer trajectory optimal control problem by using the collocation method developed in Ref.~\cite{HargravesParis1987} with a two-body gravitational model for each of the three phases: escape from the departure planet, heliocentric flight, and capture at the destination planet. 

This research is inspired by the work of Ref.~\cite{TangConway1995}.  While both Ref.~\cite{TangConway1995} and the work in this paper focus on the design of a minimum-time Earth-to-Mars interplanetary orbital transfer using low-thrust propulsion, the work presented in this paper is significantly different from the work of Ref.~\cite{TangConway1995} in the following ways. First, Ref.~\cite{TangConway1995} employs a three-phase structure consisting of Earth escape, heliocentric transfer, and Mars capture under the assumption that the planets move in circular orbits.  On the other hand, the research described in this paper employs a four-phase structure consisting of Earth and Mars alignment, Earth escape, heliocentric transfer, and Mars capture and studies the problem for the cases where the planets move in either circular or elliptic orbits (where it is noted that the alignment of the planets in the first phase of the problem ensures that the planets are in the correct position to reduce the time required to complete the interplanetary transfer).  In particular, the case of elliptic planetary motion introduces complexity into the problem due to the fact that the timing of the transfer is critical in order to obtain the minimum-time orbital transfer.  Additionally, Ref.~\cite{TangConway1995} uses two-body dynamics to model the spacecraft's motion in each phase, whereas in this study both the cases of two-body motion and motion with solar and planetary gravity perturbations are considered.  Finally, the Earth-to-heliocentric and heliocentric-to-Mars coordinate transformations employed in this paper use a different set of variables from those used in Ref.~\cite{TangConway1995} in that all angles are measured from an inertially fixed direction (specifically, this inertially fixed direction is the line of Aries).

This paper is organized as follows. Section~\ref{sect:problem-description} presents the four-phase Earth-to-Mars interplanetary orbital transfer optimal control problem by providing the modeling assumptions, phase descriptions, dynamics, path constraints, boundary conditions, variable bounds, the units used to solve the problem, the unit conversions needed to connect sequential phases, and event constraints. Section~\ref{sect:results} presents the results and discussion of the numerical optimization study using the MATLAB${}^{\textregistered}$ optimal control software, $\mathbb{GPOPS-II}$. Section~\ref{sect:conclusions} presents the conclusions on this research.  Finally, Appendix \ref{sect:coordinate-transformations} provides the derivations of the coordinate transformations used to connect the Earth escape and Mars capture phases to the heliocentric phase while Appendix \ref{sect:perturbations} provides the derivation of the gravity perturbations.

\section{Problem Description\label{sect:problem-description}}

This section develops the assumptions, constraints, boundary conditions, and equations of motion that define the four-phase Earth-to-Mars optimal control problem. Section~\ref{sect:modeling-assumptions} presents the simplifications and assumptions that were used to model the motion of the spacecraft and planets in each phase. Section \ref{sect:phase-descriptions} provides a description of the four phases of the optimal control problem.  Section~\ref{sect:equations-of-motion} covers the differential equations of motion for the spacecraft and planets as well as the path constraints that constrain the motion. Section~\ref{sect:boundary-conditions} presents the initial and terminal boundary conditions as well as the lower and upper bounds on the time, control, and state in each phase. Section~\ref{sect:scale-factors} presents the units used to solve the problem. Section~\ref{sect:unit-conversions} presents the unit conversions needed to connect sequential phases. Section~\ref{sect:event-constraints} presents the event constraints that connect the state and time of sequential phases to one another. Finally, a description of the complete optimal control problem is provided in Section~\ref{sect:optimal-control-problem}. 
 
\subsection{Modeling Assumptions\label{sect:modeling-assumptions}}

First, it is assumed that the initial true longitudes of Earth and Mars are known from ephemeris data for 1 January 2019. Second, the spacecraft and the planets move in the same inertially fixed plane. Third, the spacecraft remains in its initial parking orbit during phase 1.  As a result, phase 1 does not have any control. Fourth, in phases 2, 3, and 4 the only forces acting on the spacecraft are those due to thrust and gravitation (either two-body or multi-body gravitation, depending on the case studied). Fifth, in phases 2, 3, and 4 the control is the thrust direction. Sixth, the thrust magnitude is constant in phases 2, 3, and 4. Seventh, the mass of the spacecraft is assumed to be constant, that is, the propellant expenditure is neglected. Lastly, the longitude of the spacecraft as it starts its departure from Earth is assumed to be free.  

\subsection{Phase Descriptions\label{sect:phase-descriptions}}

Using the assumptions presented in Section~\ref{sect:modeling-assumptions}, the minimum-time trajectory optimization problem is proposed as a four-phase optimal control problem. Phase 1 is a planetary alignment phase that starts with Earth and Mars at particular locations at a specified epoch and terminates when the planets are in positions that minimize the transfer time.  Phase 2 is an Earth escape phase that begins with the spacecraft in a circular orbit relative to Earth and terminates either at the sphere of influence of Earth (for the circular and elliptic cases), or is allowed to be a free parameter in the elliptic with gravity perturbations case.  Phase 3 is a heliocentric transfer phase that begins with the spacecraft either at the sphere of influence of Earth for the circular and elliptic cases, or is allowed to be a free parameter in the elliptic with gravity perturbations case and terminates either at the sphere of influence of Mars for the circular and elliptic cases, or is allowed to be a free parameter in the elliptic with gravity perturbations case. Phase 4 is a Mars capture phase that begins either at the sphere of influence of Mars for the circular and elliptic cases, or is allowed to be a free parameter in the elliptic with gravity perturbations case and terminates with the spacecraft in a circular orbit relative to Mars.  Finally, it is noted that, for the case of elliptic planetary motion with gravity perturbations, the radii at the termination of the Earth escape phase and the start of the Mars capture phase are free because the radii at these points along the transfer can be determined as part of the optimization based on the body that has the dominant gravitational influence on the spacecraft at those points along the transfer.  

\subsection{Equations of Motion\label{sect:equations-of-motion}}

The following sets of differential equations describe the motion of Earth, Mars, and the spacecraft. First, the positions of the planets are defined by their true longitudes relative to the line of Aries. The true longitudes are governed by the following differential equations in phases $P = (1,2,3,4)$:
\begin{equation}\label{planets-equations-of-motion}
  \begin{array}{lcl}
    \dot{L}_E^{[P]}(t) & = & \displaystyle \sqrt{\frac{\mu_S^{[P]}}{\left(p_E^{[P]}\right)^3}}\left(1 + f_E\cos\left(L_E^{[P]}(t)\right) + g_E\sin\left(L_E^{[P]}(t)\right)\right)^2, \vspace{5pt}\\
    \dot{L}_M^{[P]}(t) & = & \displaystyle \sqrt{\frac{\mu_S^{[P]}}{\left(p_M^{[P]}\right)^3}}\left(1 + f_M\cos\left(L_M^{[P]}(t)\right) + g_M\sin\left(L_M^{[P]}(t)\right)\right)^2.
  \end{array}
\end{equation}
Next, the differential equations that describe the motion of the spacecraft in phases $P = (2,3,4)$ are given as:
\begin{equation}\label{spacecraft-equations-of-motion}
  \begin{array}{lcl}
    \dot{r}^{[P]}(t) & = & \displaystyle v_r^{[P]}(t), \vspace{5pt}\\
    \dot{\theta}^{[P]}(t) & = & \displaystyle \frac{v_\theta^{[P]}(t)}{r^{[P]}(t)}, \vspace{5pt}\\
    \dot{v}_r^{[P]}(t) & = & \displaystyle a^{[P]} w_r^{[P]}(t) - \frac{\mu^{[P]}}{\left(r^{[P]}(t)\right)^2} + \frac{\left(v_\theta^{[P]}(t)\right)^2}{r^{[P]}(t)} + a_{pr}^{[P]}, \vspace{5pt}\\
    \dot{v}_\theta^{[P]}(t) &  = & \displaystyle a^{[P]} w_\theta^{[P]}(t) - \frac{v_r^{[P]}(t)v_\theta^{[P]}(t)}{r^{[P]}(t)} + a_{p\theta}^{[P]}.
  \end{array}
\end{equation}
Using Cowell's method in Ref.~\cite{Bate1971}, the gravity perturbation terms are given below:
\begin{equation}\label{gravity-perturbations}
\begin{array}{lcl}
a_{pr}^{[2]} & = & \displaystyle -\mu_S^{[2]}\left[\frac{r^{[2]}}{r_{psE}^3} + \left(\frac{r_E^{[2]}}{r_{psE}^3}-\frac{1}{\left(r_E^{[2]}\right)^2}\right)\cos\left(\theta^{[2]}-L_E^{[2]}\right)\right], \vspace{5pt}\\

a_{p\theta}^{[2]} & = & \displaystyle - \mu_S^{[2]}\left(-\frac{r_E^{[2]}}{r_{psE}^3}+\frac{1}{\left(r_E^{[2]}\right)^2}\right)\sin\left(\theta^{[2]}-L_E^{[2]}\right), \vspace{5pt}\\

a_{pr}^{[3]} & = & \displaystyle -\mu_E^{[3]}\left[\frac{r^{[3]}}{r_{pe}^3} - \frac{r_E^{[3]}}{r_{pe}^3}\cos\left(\theta^{[3]}-L_E^{[3]}\right)\right] -\mu_M^{[3]}\left[\frac{r^{[3]}}{r_{pm}^3} - \frac{r_M^{[3]}}{r_{pm}^3}\cos\left(\theta^{[3]}-L_M^{[3]}\right)\right], \vspace{5pt}\\

a_{p\theta}^{[3]} & = & \displaystyle - \mu_E^{[3]}\frac{r_E^{[3]}}{r_{pe}^3}\sin\left(\theta^{[3]}-L_E^{[3]}\right) - \mu_M^{[3]}\frac{r_M^{[3]}}{r_{pm}^3}\sin\left(\theta^{[3]}-L_M^{[3]}\right), \vspace{5pt}\\

a_{pr}^{[4]} & = & \displaystyle -\mu_S^{[4]}\left[\frac{r^{[4]}}{r_{psM}^3} + \left(\frac{r_M^{[4]}}{r_{psM}^3}-\frac{1}{\left(r_M^{[4]}\right)^2}\right)\cos\left(\theta^{[4]}-L_M^{[4]}\right)\right], \vspace{5pt}\\

a_{p\theta}^{[4]} & = & \displaystyle - \mu_S^{[4]}\left(-\frac{r_M^{[4]}}{r_{psM}^3}+\frac{1}{\left(r_M^{[4]}\right)^2}\right)\sin\left(\theta^{[4]}-L_M^{[4]}\right),
\end{array} 
\end{equation}
where
\begin{equation}
\begin{array}{lcl}
r_{psE} & = & \displaystyle \sqrt{\left(r_E^{[2]}\right)^2 + \left(r^{[2]}\right)^2 + 2r^{[2]}r_E^{[2]}\cos\left(\theta^{[2]} - L_E^{[2]}\right)}, \vspace{5pt}\\

r_{pe} & = & \displaystyle \sqrt{\left(r_E^{[3]}\right)^2 + \left(r^{[3]}\right)^2 - 2r^{[3]}r_E^{[3]}\cos\left(\theta^{[3]} - L_E^{[3]}\right)}, \vspace{5pt}\\

r_{pm} & = & \displaystyle \sqrt{\left(r_M^{[3]}\right)^2 + \left(r^{[3]}\right)^2 - 2r^{[3]}r_M^{[3]}\cos\left(\theta^{[3]} - L_M^{[3]}\right)}, \vspace{5pt}\\

r_{psM} & = & \displaystyle \sqrt{\left(r_M^{[4]}\right)^2 + \left(r^{[4]}\right)^2 + 2r^{[4]}r_M^{[4]}\cos\left(\theta^{[4]} - L_M^{[4]}\right)}, \vspace{5pt}\\

\end{array}
\end{equation}
It is noted that $a_{pr}$ and $a_{p\theta}$ are both zero when solving the circular and elliptic cases because they use two-body dynamics. Table~\ref{tab:variables} summarizes the state and control variables utilized in each phase.
\begin{table}[h]
  \centering
  \renewcommand{\baselinestretch}{1.25}\footnotesize\normalfont 
  \caption{Variables used to represent the state and control in each phase. \label{tab:variables}}
  \begin{tabular}{|c|c|c|c|} \hline
    Phase & State & Control \\ \hline \hline
    Phase 1 & $\left(L_E^{[1]},L_M^{[1]}\right)$ & -- \\ \hline
    Phase 2 & $\left(r^{[2]},\theta^{[2]},v_r^{[2]},v_\theta^{[2]},L_E^{[2]},L_M^{[2]}\right)$ & $\left(w_r^{[2]},w_\theta^{[2]}\right)$ \\ \hline
    Phase 3 & $\left(r^{[3]},\theta^{[3]},v_r^{[3]},v_\theta^{[3]},L_E^{[3]},L_M^{[3]}\right)$ & $\left(w_r^{[3]},w_\theta^{[3]}\right)$ \\  \hline
    Phase 4 & $\left(r^{[4]},\theta^{[4]},v_r^{[4]},v_\theta^{[4]},L_E^{[4]},L_M^{[4]}\right)$ & $\left(w_r^{[4]},w_\theta^{[4]}\right)$ \\ \hline
  \end{tabular}
\end{table}
Next, in order to ensure that the thrust direction is a unit vector, the following equality path constraint is enforced in each of the transfer phases:
\begin{equation} \label{thrust-direction-unit-vector}
  \left(w_r^{[P]}(t)\right)^2 + \left(w_\theta^{[P]}(t)\right)^2 = 1
\end{equation}
Furthermore, the radii of the planetary orbits are given as:
 \begin{equation}\label{RadiusOfEllipse}
 \begin{array}{lcl}
  r_E^{[P]}(t) & = & \displaystyle \frac{p_E^{[P]}}{1+e_E\cos\left(\nu_E^{[P]}(t)\right)}, \vspace{3pt}\\
  r_M^{[P]}(t) & = & \displaystyle \frac{p_M^{[P]}}{1+e_M\cos\left(\nu_M^{[P]}(t)\right)},
 \end{array}
 \end{equation}
where the true anomaly of each planet is given as:
\begin{equation}\label{TrueAnomaly}
\begin{array}{lcl}
\nu_E^{[P]}(t) & = & \displaystyle L_E^{[P]}(t) - \varpi_E, \vspace{3pt}\\
\nu_M^{[P]}(t) & = & \displaystyle L_M^{[P]}(t) - \varpi_M.
\end{array}
\end{equation}
Finally, the rates of change of the radii of the planetary orbits are given as:
\begin{equation}\label{RateOfChangeRadiusOfEllipse}
 \begin{array}{lclcl}
  \dot{r}_E^{[P]}(t) & = & \displaystyle \frac{p_E^{[P]}e_E \dot{\nu}_E^{[P]}(t)\sin\left(\nu_E^{[P]}(t)\right)}{\left(1+e_E\cos\left(\nu_E^{[P]}(t)\right)\right)^2} & = & \displaystyle \frac{p_E^{[P]}e_E \dot{L}_E^{[P]}(t)\sin\left(\nu_E^{[P]}(t)\right)}{\left(1+e_E\cos\left(\nu_E^{[P]}(t)\right)\right)^2}, \vspace{3pt}\\
  \dot{r}_M^{[P]}(t) & = & \displaystyle \frac{p_M^{[P]}e_M\dot{\nu}_M^{[P]}(t)\sin\left(\nu_M^{[P]}(t)\right)}{\left(1+e_M\cos\left(\nu_M^{[P]}(t)\right)\right)^2} & = & \displaystyle \frac{p_M^{[P]}e_M \dot{L}_M^{[P]}(t)\sin\left(\nu_M^{[P]}(t)\right)}{\left(1+e_M\cos\left(\nu_M^{[P]}(t)\right)\right)^2},
 \end{array}
\end{equation}
where is it noted from Eq.~\eqref{TrueAnomaly} that
\begin{equation}\label{TrueAnomalyRate}
\begin{array}{lcl}
  \dot{\nu}_E^{[P]}(t) & = & \displaystyle\dot{L}_E^{[P]}(t), \vspace{3pt}\\
  \dot{\nu}_M^{[P]}(t) & = & \displaystyle \dot{L}_M^{[P]}(t).
\end{array}
\end{equation}
It is noted that Eq.~\eqref{RateOfChangeRadiusOfEllipse} is used as part of the event constraints that connected the various phases in the problem.  

%\clearpage

\subsection{Boundary Conditions and Bounds\label{sect:boundary-conditions}}

Bounds are placed on the time, control, and state in all four phases and are given as: 
\begin{equation}
  \begin{array}{c}
\begin{array}{lclcll}
t_{0,\min}^{[P]} & \leq & t_0^{[P]} & \leq & t_{0,\max}^{[P]}\\
t_{f,\min}^{[P]} & \leq & t_f^{[P]} & \leq & t_{f,\max}^{[P]}\\
L_{E,\min}^{[P]} & \leq & L_E^{[P]} & \leq & L_{E,\max}^{[P]}\\
L_{M,\min}^{[P]} & \leq & L_M^{[P]} & \leq & L_{M,\max}^{[P]}
  \end{array}, \quad P = (1, 2, 3, 4),\\\\
\begin{array}{lclcll}
  w_{r,\min}^{[P]} & \leq & w_r^{[P]} & \leq & w_{r,\max}^{[P]}\\
  w_{\theta,\min}^{[P]} & \leq & w_\theta^{[P]} & \leq & w_{\theta,\max}^{[P]}\\
r_{\min}^{[P]} & \leq & r^{[P]} & \leq & r_{\max}^{[P]} \\
\theta_{\min}^{[P]} & \leq & \theta^{[P]} & \leq & \theta_{\max}^{[P]} \\
v_{r,\min}^{[P]} & \leq & v_r^{[P]} & \leq & v_{r,\max}^{[P]}\\
v_{\theta,\min}^{[P]} & \leq & v_\theta^{[P]} & \leq & v_{\theta,\max}^{[P]}\\
\end{array}, \quad P = (2, 3, 4)
\end{array}
\end{equation}
All initial and terminal boundary conditions for all phases are considered to be free parameters except for the variables specified in Tables~\ref{tab:boundary_conditions_Phase_1} -- \ref{tab:boundary_conditions_Phase_4}. 
\begin{table}[h]
 \centering 
 \renewcommand{\baselinestretch}{1.25}\footnotesize\normalfont 
 \caption{Boundary conditions at start and terminus of phases 1, 2, and 4.}
 \subfloat[Boundary conditions at the start and terminus of phase 1.\label{tab:boundary_conditions_Phase_1}]{
   \phantom{--------------------------}\begin{tabular}{|c|c|}  \hline
     Variable & Value \\\hline\hline
     $t_0$ & $0$ d  \\\hline
     $L_E\left(t_0\right)$ & $101.14$~deg \\\hline
     $L_M\left(t_0\right)$ & $41.23$~deg \\\hline 
   \end{tabular}\phantom{--------------------------}
 }
 
 \subfloat[Boundary conditions at the start and terminus of phase 2\label{tab:boundary_conditions_Phase_2}]{
   \phantom{--------------------------}\begin{tabular}{|c|c|}  \hline
     Variable & Value \\\hline\hline
     $r\left(t_0\right)$ (without perturbations) & $6.6 R_E$ \\\hline
     $r\left(t_0\right)$ (with perturbations) & $6.6 R_E$ \\\hline
     $v_r\left(t_0\right)$ & 0 \\\hline
     $v_\theta\left(t_0\right)$ & $\sqrt{\mu_E/r_0}$ \\\hline
     $r\left(t_f\right)$ (without perturbations) & $R_E^{\textrm{SOI}}$ \\\hline
   \end{tabular}\phantom{--------------------------}
 }
 
 \subfloat[Boundary conditions at the start and terminus of phase 4.\label{tab:boundary_conditions_Phase_4}]{
   \phantom{--------------------------}\begin{tabular}{|c|c|}  \hline
     Variable & Value \\\hline\hline
     $r\left(t_0\right)$ (without perturbations) & $R_M^{\textrm{SOI}}$\\\hline
     $r\left(t_f\right)$ (without perturbations) & $6.0 R_M$\\\hline
     $r\left(t_f\right)$ (with perturbations) & $6.0 R_M$\\\hline
     $v_r\left(t_f\right)$ & 0 \\\hline
     $v_\theta\left(t_f\right)$ & $\sqrt{\mu_M/r_f}$\\ \hline 
   \end{tabular}\phantom{--------------------------}
 }
\end{table}
Tables~\ref{tab:bounds_Phase_1} -- \ref{tab:bounds_Phase_4} show those variables that are constrained alongside the corresponding lower and upper limits.  All other variables are free.  Finally, the physical constants and other numerical data used to model and solve the problem are given in Table~\ref{tab:physical-constants}.

\begin{table}[h]
  \caption{Lower and upper bounds during phases 1, 2, and 4.\label{tab:bounds}}
  \centering
  \renewcommand{\baselinestretch}{1.25}\footnotesize\normalfont 
  \subfloat[Lower and upper bounds for phase 1.\label{tab:bounds_Phase_1}]{
    \phantom{--------------}\begin{tabular}{|c|c|}  \hline
      Variable & [Lower Bound, Upper Bound] \\\hline\hline
      $t$ & $[0, \textrm{Free}]$  \\\hline
    \end{tabular}\phantom{--------------}}

  \subfloat[Lower and upper bounds for phase 2.\label{tab:bounds_Phase_2}]{
    \begin{tabular}{|c|c|}  \hline
      Variable & [Lower Bound, Upper Bound] \\\hline\hline
      $r$ (without perturbations) & $\left[6.6 R_E, R_E^{SOI}\right]$ \\\hline
      $r$ (with perturbations) & $[6.6 R_E, \textrm{Free}]$ \\\hline
      $v_r$ & $[0, \textrm{Free}]$ \\\hline
    \end{tabular}}

  \subfloat[Lower and upper bounds for phase 4.\label{tab:bounds_Phase_4}]{
    \begin{tabular}{|c|c|}  \hline
      Variable & [Lower Bound, Upper Bound] \\\hline\hline
      $r$ (without perturbations) & $\left[6.0 R_M, R_M^{\textrm{SOI}}\right]$ \\\hline
      $r$ (with perturbations) & $[6.0 R_M, \textrm{Free}]$ \\\hline
      $v_r$ & $[0, \textrm{Free}]$ \\\hline
    \end{tabular}}
  \end{table}

\begin{table}[h!]
  \centering
  \renewcommand{\baselinestretch}{1.25}\footnotesize\normalfont
  \caption{Physical constants.\label{tab:physical-constants}}
  \begin{tabular}{|l|c|c|} \hline
    Quantity & Value & Units \\ \hline\hline
    $R_E$ & $6.3781\times 10^{6}$ & $\textrm{m}$\\\hline
    $R_M$ & $3.3895\times 10^{6}$ & $\textrm{m}$\\\hline
    $R_E^{\textrm{SOI}}$ & $9.2455\times 10^{8}$ & $\textrm{m}$  \\\hline
    $R_M^{\textrm{SOI}}$ & $5.7717\times 10^{8}$ & $\textrm{m}$\\ \hline
    $R_{SE}$ & $1.4960 \times 10^{11}$ & $\textrm{m}$ \\\hline
    $R_{SM}$ & $2.2794 \times 10^{11}$ & $\textrm{m}$ \\\hline
    $\mu_E$   & $3.9860 \times 10^{14}$ & $\textrm{m}^3 \cdot \textrm{s}^{-2}$   \\\hline
    $\mu_S$     & $1.3271 \times 10^{20}$ & $\textrm{m}^3 \cdot \textrm{s}^{-2}$  \\\hline
   $\mu_M$   & $4.2828 \times 10^{13}$ & $\textrm{m}^3 \cdot \textrm{s}^{-2}$  \\\hline
   $\varpi_E$ & 102.9 & deg \\\hline
   $\varpi_M$ & 336.0 & deg \\\hline
  \end{tabular}
\end{table}

% \clearpage

\subsection{Scale Factors\label{sect:scale-factors}}

The units used in each phase $P\in (1,2,3,4)$ of the four-phase optimal control problem were chosen such that the gravitational parameter of the central body in that phase is equal to unity.   In order to attain a gravitational parameter of unity in each phase, the distance unit, $D^{[P]},\; (P=1,2,3,4)$, speed unit, $V^{[P]},\; (P=1,2,3,4)$, and time unit, $T^{[P]},\; (P=1,2,3,4)$, were chosen as follows:
\begin{equation}
  \left[
    \begin{array}{c}
      D^{[P]} \\
      V^{[P]} \\
      T^{[P]}
    \end{array}
  \right]
  =
  \left[
    \begin{array}{c}
      R^{[P]} \\
      \sqrt{\mu^{[P]}/D^{[P]}} \\
      \sqrt{\left(D^{[P]}\right)^3/\mu^{[P]}}
    \end{array}
  \right], \quad P\in (1,2,3,4),
\end{equation}
where
\begin{equation}
  \left[
    \begin{array}{c}
      R^{[1]} \\
      R^{[2]} \\
      R^{[3]} \\
      R^{[4]}
    \end{array}
  \right]
  =
  \left[
    \begin{array}{c}
      R_{SE} \\
      R_E \\
      R_{SE} \\
      R_M
    \end{array}
  \right].
\end{equation}

% \clearpage

\subsection{Unit Conversions\label{sect:unit-conversions}}

A set of unit conversions are created to connect all variables across the sequential phases with identical units provided in Section~\ref{sect:scale-factors}. First, in the Alignment-to-Earth event constraint, the variables at the start of phase 2 are converted from Earth units $\left(D^{[2]} , V^{[2]} , T^{[2]}\right)$ to heliocentric units $\left(D^{[1]} , V^{[1]} , T^{[1]}\right)$ in order to connect the terminus of phase 1 to the start of phase 2. Next, in the Earth-to-Heliocentric event constraint, the variables at the terminus of phase 2 are converted from Earth units $\left(D^{[2]} , V^{[2]} , T^{[2]}\right)$ to heliocentric units $\left(D^{[3]} , V^{[3]} , T^{[3]}\right)$ in order to connect the terminus of phase 2 to the start of phase 3. Lastly, in the Heliocentric-to-Mars event constraint, the variables at the start of phase 4 are converted from Mars units $\left(D^{[4]} , V^{[4]} , T^{[4]}\right)$ to heliocentric units $\left(D^{[3]} , V^{[3]} , T^{[3]}\right)$ in order to connect the terminus of phase 3 to the start of phase 4. The necessary conversion factors to transform Earth units to heliocentric units are given as:

\begin{equation}\label{e2hunitconversions}
  \left[
    \begin{array}{c}
      D_{SE} \vspace{5pt}\\
      V_{SE} \vspace{5pt}\\
      T_{SE}
    \end{array}
  \right]
  =
  \left[
    \begin{array}{c}
      D^{[2]}/D^{[1]} \vspace{5pt}\\
      V^{[2]}/V^{[1]} \vspace{5pt}\\
      T^{[2]}/T^{[1]}
    \end{array}
  \right]
  =
  \left[
    \begin{array}{c}
      D^{[2]}/D^{[3]} \vspace{5pt}\\
      V^{[2]}/V^{[3]} \vspace{5pt}\\
      T^{[2]}/T^{[3]}
    \end{array}
  \right].
\end{equation}
The necessary conversion factors to transform Mars units to heliocentric units are given as:
\begin{equation}\label{m2hunitconversions}
  \left[
    \begin{array}{c}
      D_{SM} \vspace{5pt}\\
      V_{SM} \vspace{5pt}\\
      T_{SM}
    \end{array}
  \right]
  =
  \left[
    \begin{array}{c}
      D^{[4]}/D^{[1]} \vspace{5pt}\\
      V^{[4]}/V^{[1]} \vspace{5pt}\\
      T^{[4]}/T^{[1]} 
    \end{array}
  \right]
  =
  \left[
    \begin{array}{c}
      D^{[4]}/D^{[3]} \vspace{5pt}\\
      V^{[4]}/V^{[3]} \vspace{5pt}\\
      T^{[4]}/T^{[3]}
    \end{array}
  \right].
\end{equation}

% \clearpage

\subsection{Event Constraints\label{sect:event-constraints}}

\begin{figure}[h!]
  \centering
   \includegraphics[width=1.75in]{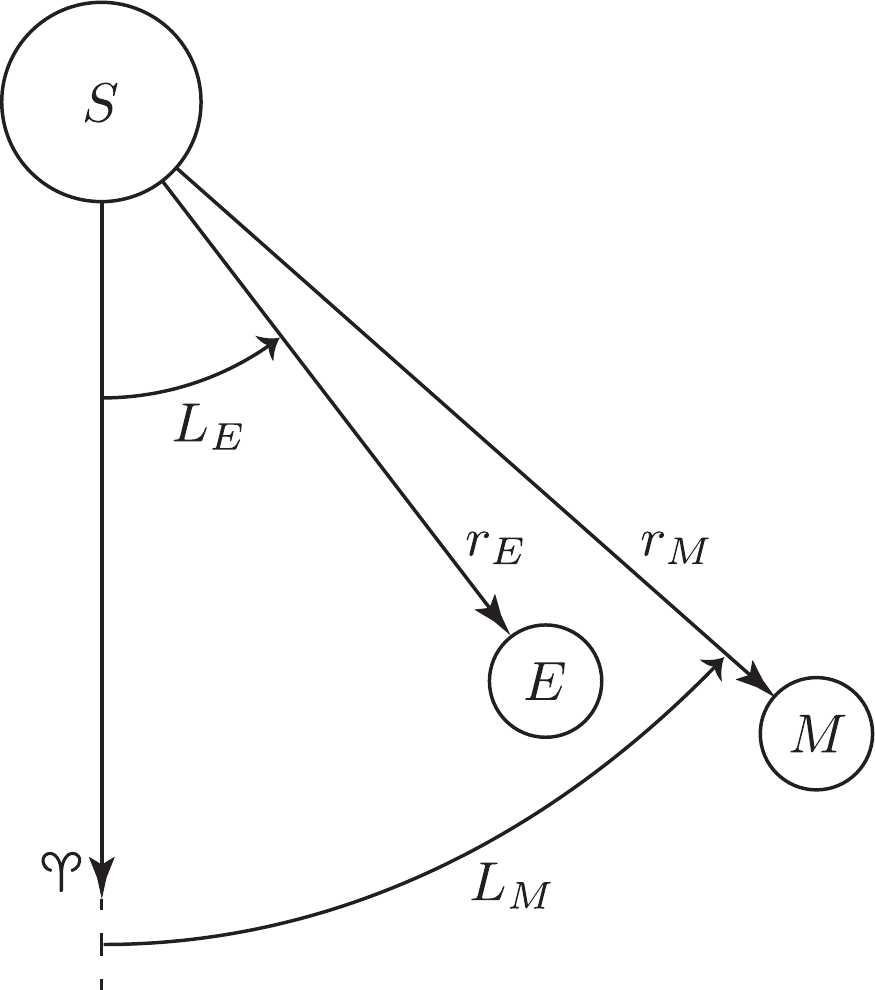}
  \caption{Geometry of the planetary alignment and Earth escape phases.}
  \label{fig:planets-earth-geometry}
\end{figure}

\begin{figure}[h!]
  \centering
  \includegraphics[width=3.25in]{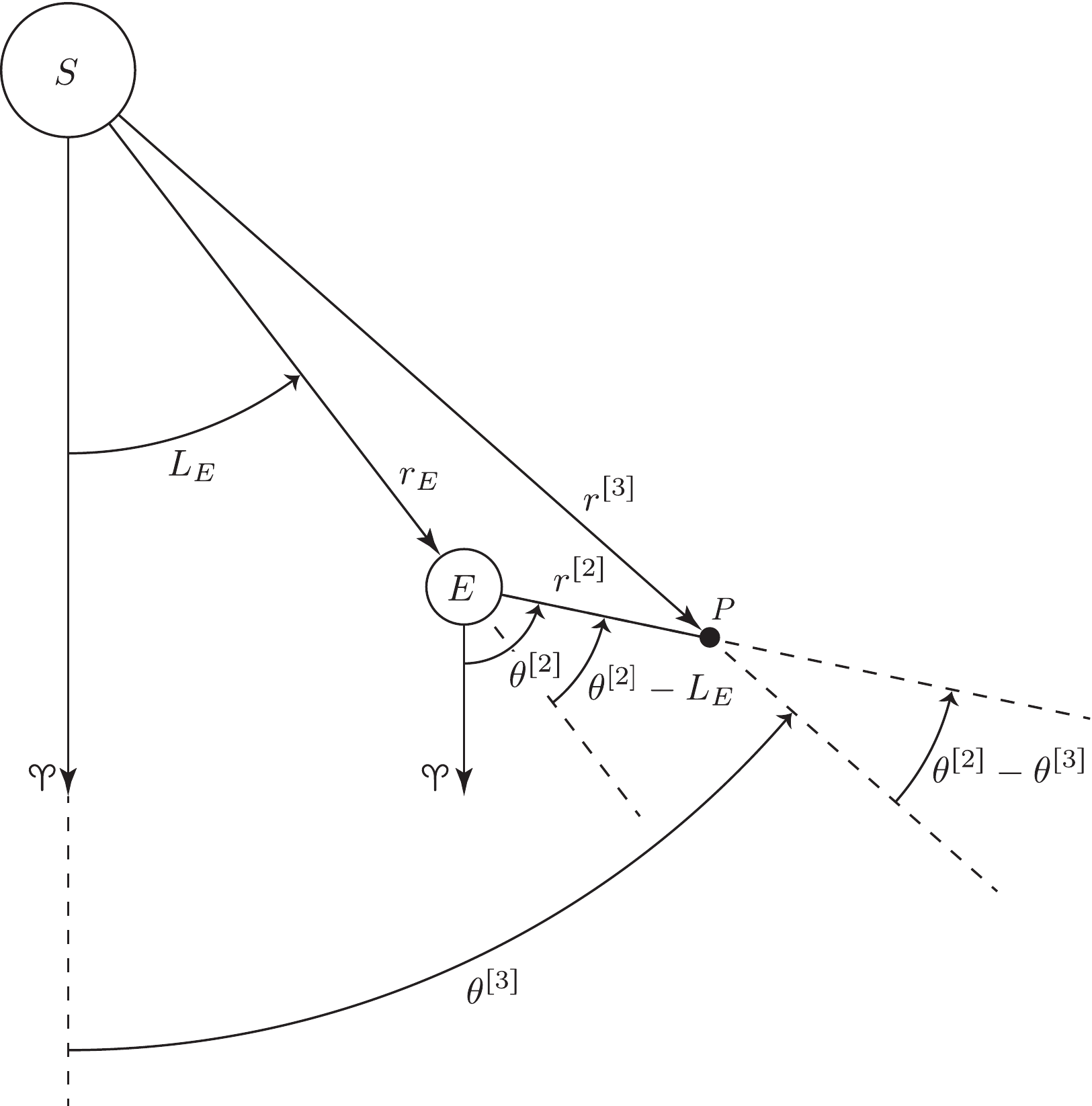}
  \caption{Geometry of the Earth escape and heliocentric transfer phases.}
  \label{fig:earth-heliocentric-geometry}
\end{figure}

\begin{figure}[h!]
  \centering
    \includegraphics[width=3.25in]{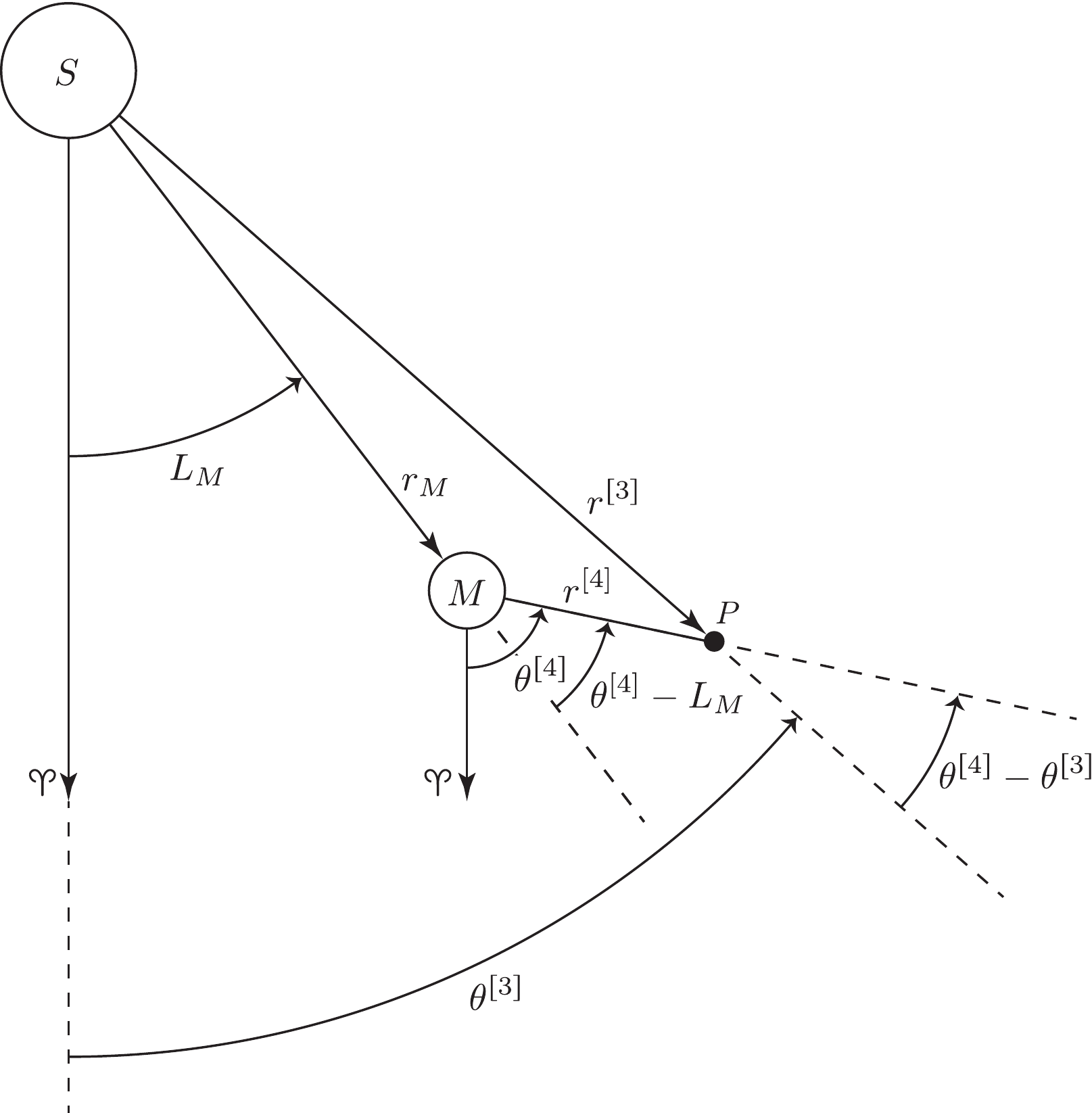}
  \caption{Geometry of the heliocentric transfer and Mars capture phases.}
  \label{fig:heliocentric-mars-geometry}
\end{figure}

Event constraints are utilized to connect subsequent phases to one another by enforcing the constraints at the endpoints, so that the variables at the terminus of one phase can be transformed to the variables at the start of the subsequent phase. The event constraints in this research connect phase 1 to phase 2, phase 2 to phase 3, and phase 3 to phase 4. The following variables are utilized in the event constraints to ensure continuity: time, true longitude of Earth and Mars, and the spacecraft's position and velocity components. The unit conversions stated in Section~\ref{sect:unit-conversions} are implemented in the event constraints to provide continuous units. It is noted in this research that the event constraints are formulated different from those found in Ref.~\cite{TangConway1995} and is constructed in such a way to allow for non-circular planetary motion unlike in Ref.~\cite{TangConway1995}.  First, the Alignment-to-Earth event constraints that connect the terminus of phase 1 to the start of phase 2 are given as:
\begin{equation}\label{events-planetary-to-earth}
  \begin{array}{lcl}
    t_f^{[1]} & = & t_0^{[2]}, \vspace{3pt}\\
    L_E^{[1]}\left(t_f^{[1]}\right) & = & L_E^{[2]}\left(t_0^{[2]}\right), \vspace{3pt}\\
    L_M^{[1]}\left(t_f^{[1]}\right) & = & L_M^{[2]}\left(t_0^{[2]}\right).
  \end{array}
\end{equation}
The geometry of phase 1 and phase 2 is shown in Fig.~\ref{fig:planets-earth-geometry}.
Second, the event constraints that connect the terminus of phase 2 to the start of phase 3 are given as (see \ref{sect:coordinate-transformations} for derivation of the coordinate transformations):
\begin{equation}\label{events-earth-to-heliocentric}
  \begin{array}{lcl}
    t_f^{[2]}T_{SE} & = & t_0^{[3]}, \vspace{3pt} \\
    
    L_E^{[2]}\left(t_f^{[2]}\right) & = & L_E^{[3]}\left(t_0^{[3]}\right), \vspace{3pt}\\
    
    L_M^{[2]}\left(t_f^{[2]}\right) & = & L_M^{[3]}\left(t_0^{[3]}\right), \vspace{3pt}\\
    
    r^{[2]}\left(t_f^{[2]}\right)D_{SE} & = & -r_E^{[3]}\left(t_0^{[3]}\right) \cos(\alpha) + r^{[3]}\left(t_0^{[3]}\right)\cos(\beta), \vspace{3pt}\\
    
    0 & = & -r_E^{[3]}\left(t_0^{[3]}\right) \sin(\alpha) + r^{[3]}\left(t_0^{[3]}\right) \sin(\beta), \vspace{3pt}\\
    
    v_r^{[2]}\left(t_f^{[2]}\right)V_{SE} & = & -\dot{r}_E^{[3]}\left(t_0^{[3]}\right)\cos(\alpha)-r_E^{[3]}\left(t_0^{[3]}\right) \dot{L}_E^{[3]}\left(t_0^{[3]}\right) \sin(\alpha) \\
    & & + v_r^{[3]}\left(t_0^{[3]}\right)\cos(\beta) + v_\theta^{[3]}\left(t_0^{[3]}\right) \sin(\beta), \vspace{3pt}\\
     
    v_\theta^{[2]}\left(t_f^{[2]}\right)V_{SE} & = & -\dot{r}_E^{[3]}\left(t_0^{[3]}\right)\sin(\alpha)-r_E^{[3]}\left(t_0^{[3]}\right) \dot{L}_E^{[3]}\left(t_0^{[3]}\right) \cos(\alpha) \\
    & & - v_r^{[3]}\left(t_0^{[3]}\right)\sin(\beta) + v_\theta^{[3]}\left(t_0^{[3]}\right) \cos(\beta),
    
  \end{array}
\end{equation}
where
\begin{equation}
\begin{array}{lcl}
\alpha & = & \theta^{[2]}\left(t_f^{[2]}\right)-L_E^{[3]}\left(t_0^{[3]}\right), \vspace{3pt}\\
\beta & = & \theta^{[2]}\left(t_f^{[2]}\right)-\theta^{[3]}\left(t_0^{[3]}\right).
\end{array}
\end{equation}
The geometry of phase 2 and phase 3 is shown in Fig.~\ref{fig:earth-heliocentric-geometry}.
Finally, the event constraints that connect the terminus of phase 3 to the start of phase 4 are given as (see \ref{sect:coordinate-transformations} for the derivation of the coordinate transformations):
\begin{equation}\label{events-heliocentric-to-mars}
  \begin{array}{lcl}
    t_0^{[4]}T_{SM} & = & t_f^{[3]}, \vspace{3pt}\\
    
    L_E^{[4]}\left(t_0^{[4]}\right) & = & L_E^{[3]}\left(t_f^{[3]}\right), \vspace{3pt}\\
    
    L_M^{[4]}\left(t_0^{[4]}\right) & = & L_M^{[3]}\left(t_f^{[3]}\right), \vspace{3pt}\\
    
    r^{[4]}\left(t_0^{[4]}\right)D_{SM} & = & -r_M^{[3]}\left(t_f^{[3]}\right) \cos(\gamma) + r^{[3]}\left(t_f^{[3]}\right)\cos(\delta), \vspace{3pt}\\
    
    0 & = & -r_M^{[3]}\left(t_f^{[3]}\right) \sin(\gamma) + r^{[3]}\left(t_f^{[3]}\right) \sin(\delta), \vspace{3pt}\\
    
    v_r^{[4]}\left(t_0^{[4]}\right)V_{SM} & = & -\dot{r}_M^{[3]}\left(t_f^{[3]}\right)\cos(\gamma)-r_M^{[3]}\left(t_f^{[3]}\right) \dot{L}_M^{[3]}\left(t_f^{[3]}\right) \sin(\gamma) \\
    & & + v_r^{[3]}\left(t_f^{[3]}\right)\cos(\delta) + v_\theta^{[3]}\left(t_f^{[3]}\right) \sin(\delta), \vspace{3pt}\\
    
    v_\theta^{[4]}\left(t_0^{[4]}\right)V_{SM} & = & \dot{r}_M^{[3]}\left(t_f^{[3]}\right)\sin(\gamma)-r_M^{[3]}\left(t_f^{[3]}\right) \dot{L}_M^{[3]}\left(t_f^{[3]}\right) \cos(\gamma) \\
    & & - v_r^{[3]}\left(t_f^{[3]}\right)\sin(\delta) + v_\theta^{[3]}\left(t_f^{[3]}\right) \cos(\delta),
    
  \end{array}
\end{equation}
\begin{equation}
\begin{array}{lcl}
\gamma & = & \theta^{[4]}\left(t_0^{[4]}\right)-L_M^{[3]}\left(t_f^{[3]}\right), \vspace{3pt}\\
\delta & = & \theta^{[4]}\left(t_0^{[4]}\right)-\theta^{[3]}\left(t_f^{[3]}\right).
\end{array}
\end{equation}
The geometry of phase 3 and phase 4 is shown in Fig.~\ref{fig:heliocentric-mars-geometry}.

%\clearpage

%
%

\subsection{Optimal Control Problem\label{sect:optimal-control-problem}}

For the Earth-to-Mars transfer, the optimal control problem is stated as follows. Determine the state $\left(L_E^{[P]}, L_M^{[P]}\right)$,~$P = 1$ and $\left(r^{[P]},\theta^{[P]},v_r^{[P]},v_\theta^{[P]},L_E^{[P]},L_M^{[P]}\right)$,~$P = (2,3,4)$, the control $\left(w_r^{[P]},w_\theta^{[P]}\right)$,~$P = (2,3,4)$, as well as the initial and terminal times $\left(t_0^{[P]},t_f^{[P]}\right)$,~$P = (1,2,3,4)$ which minimize
\begin{equation}\label{objective-functional}
J = t_f^{[4]} - t_0^{[2]}
\end{equation}
while satisfying the dynamics and path constraints in Section~\ref{sect:equations-of-motion}, the variable bounds and boundary conditions in Section~\ref{sect:boundary-conditions}, and the event constraints in Section~\ref{sect:event-constraints}.

\section{Results and Discussion\label{sect:results}}

\begin{figure}[h]
  \centering
  \includegraphics[width=3.25in]{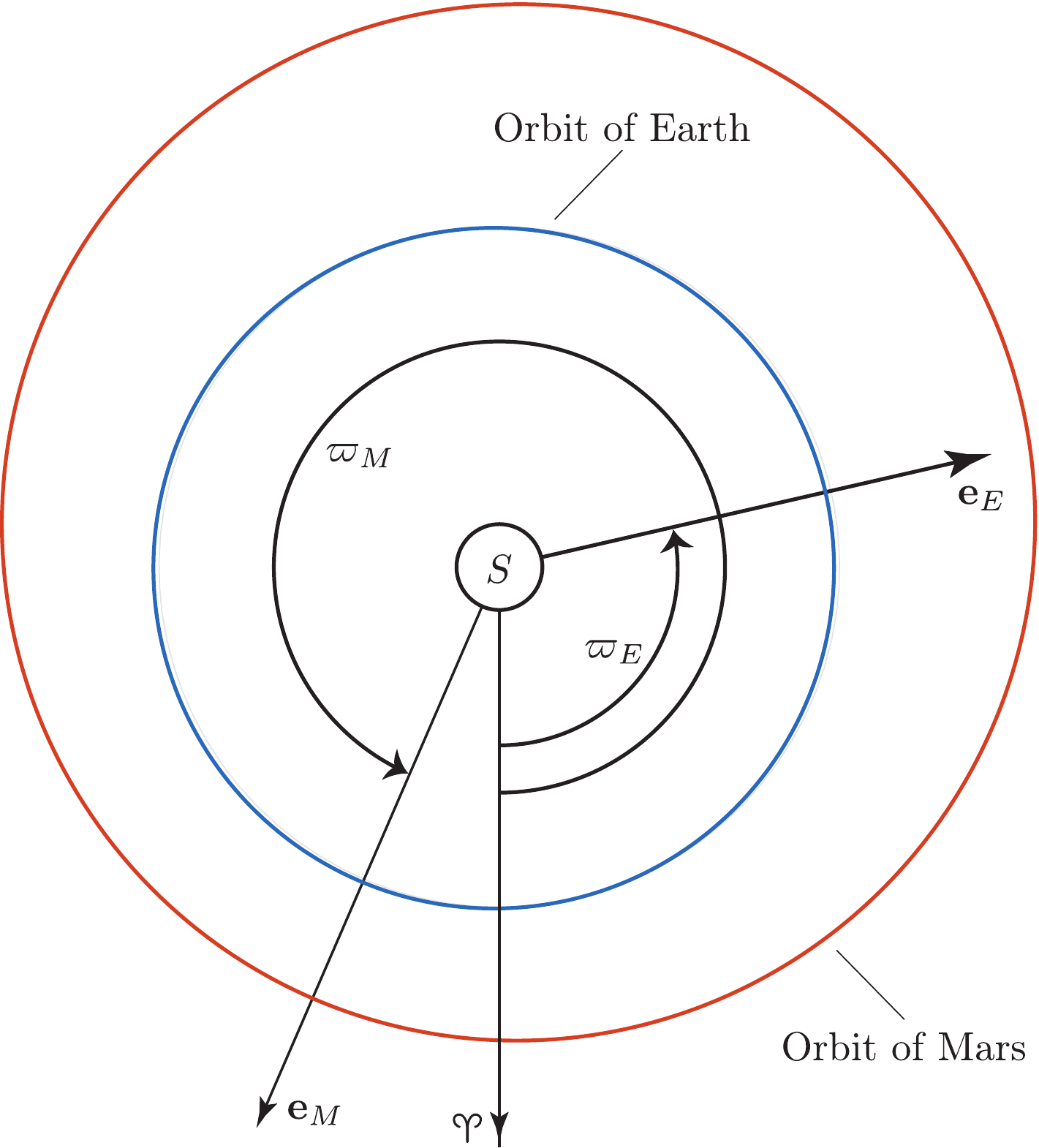}
  \caption{Relative location of Earth and Mars on 1 January 2019 at 00:00 UTC.\label{fig:planet-ellipses}}
\end{figure}

This section presents the results acquired by solving the four-phase Earth-to-Mars orbital transfer optimal control problem described in Section~\ref{sect:problem-description}.  For all results presented in this research, it is assumed that the planetary alignment phase starts with Earth and Mars at the locations defined by the NASA HORIZONS J2000 basis \cite{NASAHORIZONS} on 1 January 2019 at 00:00 Coordinated Universal Time (UTC) and that the longitudes of perihelion of Earth and Mars, $\varpi_E$ and $\varpi_M$, are $102.9$~deg and $336.0$~deg, respectively (see Fig.~\ref{fig:planet-ellipses}).  Next, the cases of circular planetary motion, elliptic planetary motion, and elliptic planetary motion with gravity perturbations are considered.  For the case of elliptic planetary motion, it is assumed that eccentricities of Earth and Mars orbit are $e_E = 0.0167$ and $e_M = 0.0935$.  Solutions are obtained for the three different cases of planetary motion for $a=(9.8,9.9,10.0,10.1,10.2)\times 10^{-4}~\textrm{m}\cdot\textrm{s}^{-2}$.  Finally, the particular case of elliptic planetary motion with gravity perturbations and $a = 9.8\times10^{-4}~\textrm{m}\cdot\textrm{s}^{-2}$ is used to show the key features of the optimized solutions.

The Earth-to-Mars transfer problem is solved using the general-purpose MATLAB${}^{\textregistered}$ optimal control software $\mathbb{GPOPS-II}$ \cite{PattersonRao2014} with the following settings.  First, the mesh refinement method described in Ref.~\cite{LiuHagerRao2018} was used with a mesh refinement accuracy tolerance $10^{-7}$.  Second the nonlinear programming problem (NLP) solver {\em IPOPT} \cite{BieglerZavala2009} was employed in full Newton (second derivative) mode with an NLP solver tolerance $\epsilon_{\textrm{NLP}}$, of $10^{-10}$, and all first and second derivatives required by {\em IPOPT} were obtained using the open-source algorithmic differentiations software {\em ADiGator}\cite{WeinsteinRao2017}.  All computations were performed using a 2.9 GHz Intel Core i9 MacBook Pro running Mac OS version 10.13.6 (High Sierra) with 32GB 2400MHz DDR4 RAM and MATLAB${}^{\textregistered}$ Version R2018b (build 9.5.0.944444).

Next, an initial guess must be supplied in order to solve the four-phase optimal control problem using $\mathbb{GPOPS-II}$.  In this research the following initial guesses were generated using $\mathbb{GPOPS-II}$ for each of the four phases of the problem.  First, for the planetary alignment phase $\mathbb{GPOPS-II}$ was used to determine the planetary trajectories starting at $(L_{E}(0), L_{M}(0))=(101.14,41.23)$~deg and terminating $150$~d later.  Second, for the Earth escape phase $\mathbb{GPOPS-II}$ was used to determine the minimum-time trajectory and control starting with the spacecraft in an Earth circular orbit of radius $6.6 R_E$ and terminating in a circular orbit of radius $R_E^{SOI}$ relative to the Earth.  Third, for the heliocentric phase $\mathbb{GPOPS-II}$ was used to determine the minimum-time trajectory and control starting in a heliocentric circular orbit of radius $1$~AU and terminating with the spacecraft in a heliocentric circular orbit with a radius $1.5$~AU.  Finally, for the Mars capture phase $\mathbb{GPOPS-II}$ was used to determine the trajectory and control starting with the spacecraft in Mars circular orbit of radius $R_M^{SOI}$ and terminating in Mars circular orbit of radius $6.0 R_M$.  Each of the aforementioned four one-phase solutions were then used as an initial guess for the four-phase circular planetary motion problem with $a=9.8\times 10^{-4}~\textrm{m}\cdot\textrm{s}^{-2}$.  Finally, the solution obtained for each case of planetary motion and each value of $a$ was obtained using the previously obtained solution as an initial guess.  

\subsection{Comparison with Results Obtained in Ref.~\cite{TangConway1995}\label{sect:TangConwayComparison}}

This section provides a comparison of the results obtained in Ref.~\cite{TangConway1995} against the formulation developed in this paper.  In order to provide a direct comparison with the results obtained in Ref.~\cite{TangConway1995}, the following modifications are made to the problem formulation given in Section \ref{sect:problem-description} of this paper.  First, because Ref.~\cite{TangConway1995} does not include a planetary alignment phase, the problem formulation given Section \ref{sect:problem-description} is simplified to exclude the planetary alignment phase (that is, only Earth escape, heliocentric, and Mars capture phases are included).  Second, it is assumed that the planetary motion is circular and that gravity perturbations are excluded.   Third, the values of thrust specific force are $a=(9.604,9.702,9.800,10.290,10.780)\times 10^{-4}~\textrm{m}\cdot\textrm{s}^{-2}$ (and are those used in Ref.~\cite{TangConway1995}).  Finally, the initial planetary phase angles are fixed such that Mars leads Earth by $0.9666~\textrm{rad}\approx 55~\textrm{deg}$ (where $0.966~\textrm{rad}$ is the value obtained in Ref.~\cite{TangConway1995} using local optimization and is the value for which results in Ref.~\cite{TangConway1995} are shown).  Table~\ref{tab:TC-circ-acceleration-times} shows the transfer times obtained in Ref.~\cite{TangConway1995} against the transfer times obtained using the modified problem formulation developed in this paper.  It is seen that the transfer times obtained in this research are smaller than the transfers times obtained in Ref.~\cite{TangConway1995} where the difference in transfer time decreases as $a$ increases.  

\begin{table}  
\centering
\renewcommand{\baselinestretch}{1.25}\footnotesize\normalfont
\caption{Transfer times obtained from Ref.~\cite{TangConway1995} and our method.\label{tab:TC-circ-acceleration-times}}
\centering
  \begin{tabular}{|c|c|c|} \hline
    $a~(\textrm{m}\cdot \textrm{s}{}^{-2})\times 10^{-4}$ & Ref.~\cite{TangConway1995} Transfer Time (d) & Transfer Time (d)  \\\hline\hline 
    9.604 & 229.15 & 223.60 \\\hline
    9.702 & 224.37 & 222.74 \\\hline
    9.800 & 222.14 & 221.89 \\\hline
    10.290 & 218.18 & 217.84 \\\hline
    10.780 & 214.51 & 214.07 \\ \hline
  \end{tabular}
\end{table}

\subsection{Overall Performance\label{sect:overall-performance}}

While Section \ref{sect:TangConwayComparison} provides a comparison with the work of Ref.~\cite{TangConway1995}, the remainder of this study focuses on the performance obtained using the problem formulation developed in this research.  Specifically, the minimum-time results obtained in this study for the circular, elliptic, and elliptic with gravity perturbations cases are shown in Tables~\ref{tab:circ-acceleration-times}~--~\ref{tab:pert-acceleration-times}, respectively for $a=(9.8,9.9,10.0,10.1,10.2)\times 10^{-4}~\textrm{m}\cdot\textrm{s}^{-2}$ (where it is noted that the values of $a$ used in this study differ from those used in Ref.~\cite{TangConway1995}).  For each value of $a$, the optimal duration for each phase of the problem along with the total minimum interplanetary transfer time and the corresponding start date for the transfer are presented. It is seen from Tables~\ref{tab:circ-acceleration-times} and~\ref{tab:ell-acceleration-times} that, as $a$ increases, the time for the planetary alignment phase also increases while the time for each of the transfer phases and the total transfer time decreases. Furthermore, Table~\ref{tab:pert-acceleration-times} shows that, as $a$ increases, the time for the planetary alignment phase also increases, the transfer phases do not follow a distinct pattern for either increasing or decreasing, and the total transfer time decreases. It is noted that, when gravity perturbations are included, the phases do not follow a pattern similar to those of the circular and elliptic cases because when gravity perturbations are included the radius at the terminus of phase 2 and start of phase 4 are free.

Next, when comparing the results of the circular and elliptic cases obtained in this study, it is seen in Tables~\ref{tab:circ-acceleration-times} and~\ref{tab:ell-acceleration-times} that, for all values of $a$, the elliptic case has a shorter phase 1 duration by an average of $1.63~\textrm{d}$, a longer phase 2 duration by approximately $0.01~\textrm{d}$, a shorter phase 3 duration by an average of $18.98~\textrm{d}$, a shorter phase 4 duration by an average of $0.01~\textrm{d}$, and a shorter total transfer time by an average of $18.98~\textrm{d}$.  When comparing the results of the elliptic and elliptic with gravity perturbations cases obtained in this study, Tables~\ref{tab:ell-acceleration-times} and~\ref{tab:pert-acceleration-times}  show that, for all values of $a$, the elliptic with gravity perturbations case has a shorter phase 1 duration by an average of $1.13~\textrm{d}$.  Furthermore, because the transfer phases of the elliptic with gravity perturbations case do not follow an increasing or decreasing pattern, some of the phases have a longer or shorter duration when compared with the elliptic case depending upon the value of $a$.  Consequently, the phase 2 duration was found to differ by an average of $0.83~\textrm{d}$, the phase 3 duration differed by an average of $1.33~\textrm{d}$, and the phase 4 duration differed by an average of $0.32~\textrm{d}$. For all values of $a$, the elliptic with gravity perturbations case has a longer total transfer time by an average of $1.83~\textrm{d}$.  Finally, as $a$ increases, the start date of the transfer for the circular case lies between 1 July 2020 to 3 July 2020, the elliptic case lies between 30 June 2020 to 2 July 2020, and the elliptic with gravity perturbations case lies between 28 June 2020 to 1 July 2020.

\begin{table}  
\centering
\renewcommand{\baselinestretch}{1.25}\footnotesize\normalfont
\caption{Transfer times obtained in this research for the cases of circular planetary motion, elliptic planetary motion, and elliptic planetary motion with gravity perturbations.\label{tab:transfer-times}}
\centering

\subfloat[Transfer times for circular planetary motion.\label{tab:circ-acceleration-times}]{
    \begin{tabular}{|c|c|c|c|c|c|c|} \hline
    $a~(\textrm{m}\cdot \textrm{s}{}^{-2})\times 10^{-4}$& Phase 1 (d) & Phase 2 (d) & Phase 3 (d) & Phase 4 (d) & Transfer Time (d) & Start Date \\\hline\hline 
   9.8 & 547.63 & 33.27 & 162.48 & 19.31 & 215.05 & 1 July 2020 \\\hline
    9.9 & 548.20 & 32.98 & 161.76 & 19.16 & 213.90 & 2 July 2020 \\\hline
    10.0 & 548.75 & 32.69 & 161.06 & 19.01 & 212.76 & 2 July 2020 \\\hline
    10.1 & 549.30 & 32.40 & 160.37 & 18.87 & 211.64 & 3 July 2020 \\\hline
    10.2 & 549.84 & 32.11 & 159.70 & 18.73 & 210.53 & 3 July 2020 \\ \hline
  \end{tabular}}

\subfloat[Transfer times for elliptic planetary motion. \label{tab:ell-acceleration-times}]{
  \begin{tabular}{|c|c|c|c|c|c|c|} \hline
    $a~(\textrm{m}\cdot \textrm{s}{}^{-2})\times 10^{-4}$& Phase 1 (d) & Phase 2 (d) & Phase 3 (d) & Phase 4 (d) & Transfer Time (d) & Start Date \\\hline\hline 
    9.8 & 546.04 & 33.28 & 143.36 & 19.30 & 195.94 & 30 June 2020 \\\hline
    9.9 & 546.58 & 32.99 & 142.72 & 19.15 & 194.86 & 30 June 2020 \\\hline
    10.0 & 547.12 & 32.70 & 142.08 & 19.00 & 193.78 & 1 July 2020 \\\hline
    10.1 & 547.66 & 32.41 & 141.46 & 18.86 & 192.72 & 1 July 2020 \\\hline
    10.2 & 548.19 & 32.12 & 140.84 & 18.72 & 191.68 & 2 July 2020 \\ \hline
  \end{tabular}
}

\subfloat[Transfer times for elliptic planetary motion with gravity perturbations.\label{tab:pert-acceleration-times}]{
  \begin{tabular}{|c|c|c|c|c|c|c|} \hline
    $a~(\textrm{m}\cdot \textrm{s}{}^{-2})\times 10^{-4}$& Phase 1 (d) & Phase 2 (d) & Phase 3 (d) & Phase 4 (d) & Transfer Time (d) & Start Date \\\hline\hline 
    9.8  & 544.88 & 33.31 & 146.25 & 18.27 & 197.83 & 28 June 2020 \\\hline
    9.9  & 545.44 & 33.05 & 146.27 & 17.39 & 196.71 & 29 June 2020 \\\hline
    10.0 & 546.00 & 36.76 & 138.97 & 19.89 & 195.61 & 30 June 2020 \\\hline
    10.1 & 546.55 & 32.43 & 144.31 & 17.79 & 194.53 & 30 June 2020 \\\hline
    10.2 & 547.09 & 32.08 & 141.29 & 20.09 & 193.46 & 1 July 2020 \\ \hline
  \end{tabular}
}

\end{table}

\subsection{Key Features of Optimized Solutions\label{sect:key-features}}

This section shows the key features of all optimized solutions using the particular case $a=9.8\times 10^{-4}~\textrm{m}\cdot\textrm{s}^{-2}$ for the case of elliptic planetary motion with gravity perturbations. These key features of the solutions are shown for each phase of the four-phase problem.  Particular attention is given to the behavior of the spacecraft trajectory and the control that produces that behavior. 

It was found from the NASA HORIZONS J2000 database \cite{NASAHORIZONS} that the values of the true longitudes of Earth and Mars were $101.14$~deg and $41.23$~deg, respectively. Next, during the planetary alignment phase the Earth traverses approximately $1.50$ orbits about the Sun while Mars traverses approximately $0.79$ orbits about the Sun.  In addition, during this phase the spacecraft remains in its initial orbit (that is, no propulsive force is exerted on the spacecraft). The optimized duration of the planetary alignment phase is $544.88~\textrm{d}$ and terminates on 28 June 2020.

\begin{figure}[h]
  \centering
  \subfloat[Optimal two-dimensional trajectory.\label{fig:pert-Trajectory_Phase_2}]{\includegraphics[width=3.25in]{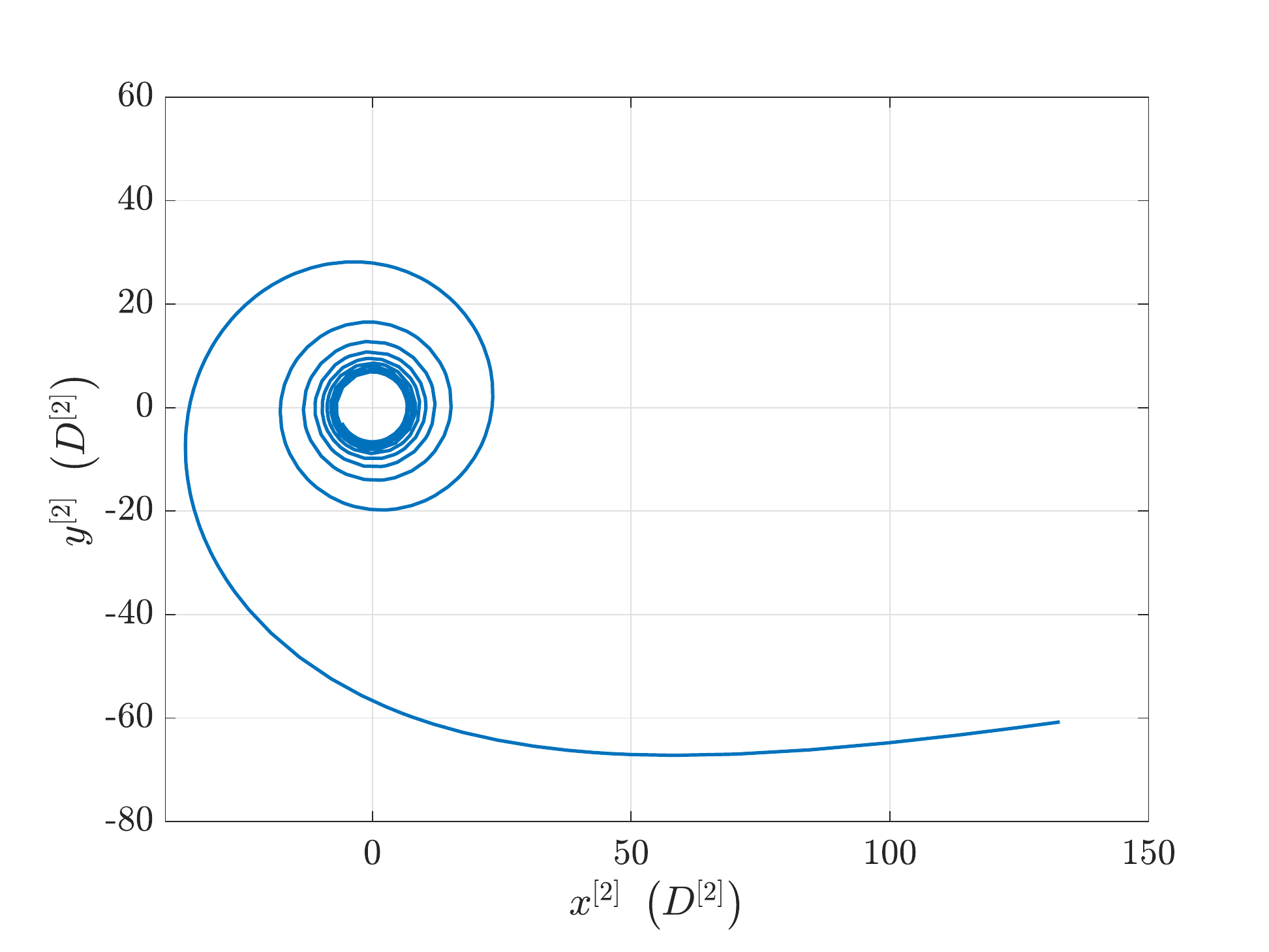}}~~\subfloat[Spacecraft eccentricity. \label{fig:pert-eSC_Phase_2}]{\includegraphics[width=3.25in]{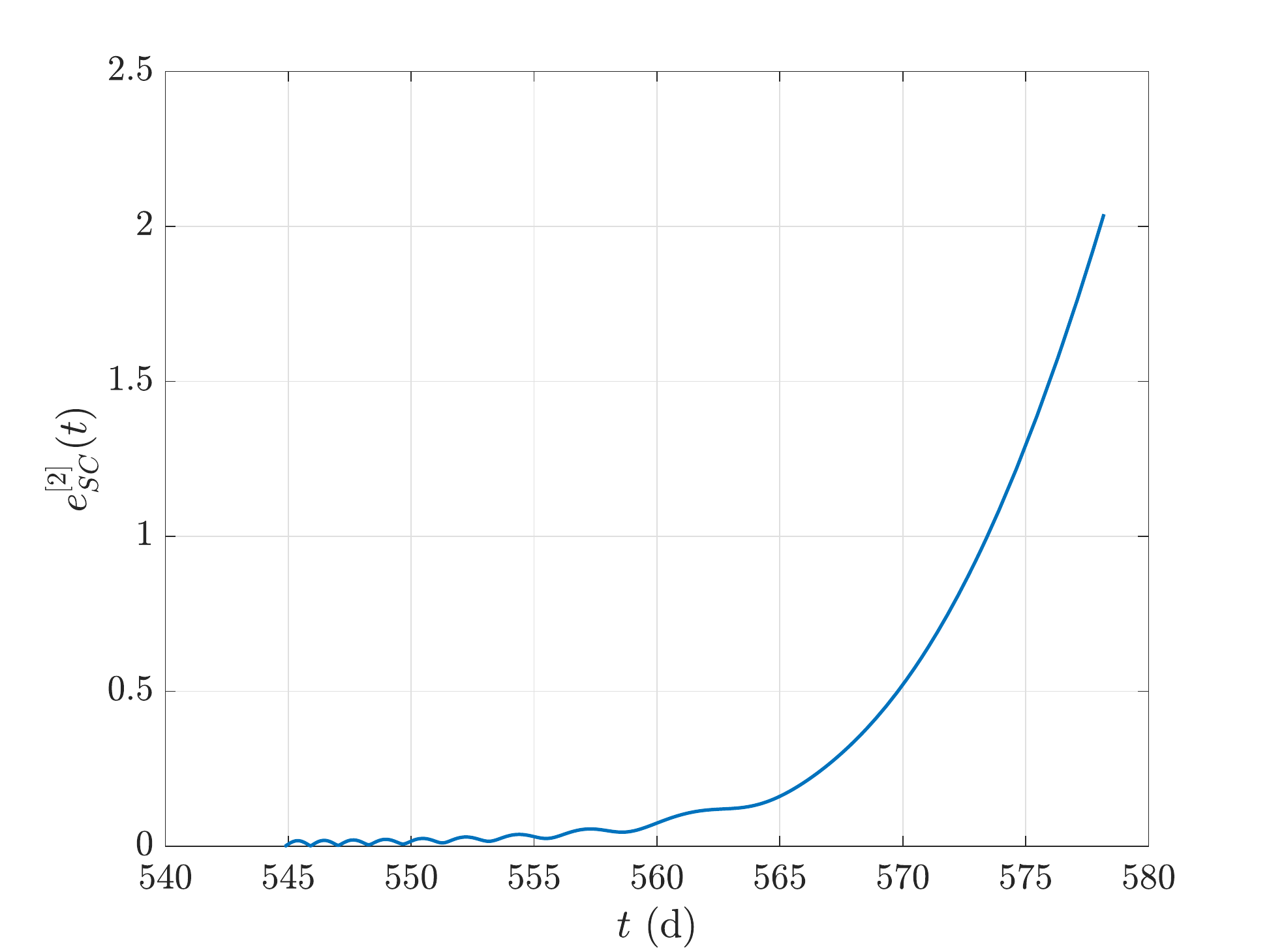}}

  \subfloat[Optimal controls $\left(w_r^{[2]},w_\theta^{[2]}\right)$ vs. $t$. \label{fig:pert-Control_Components_Phase_2}]{\includegraphics[width=3.25in]{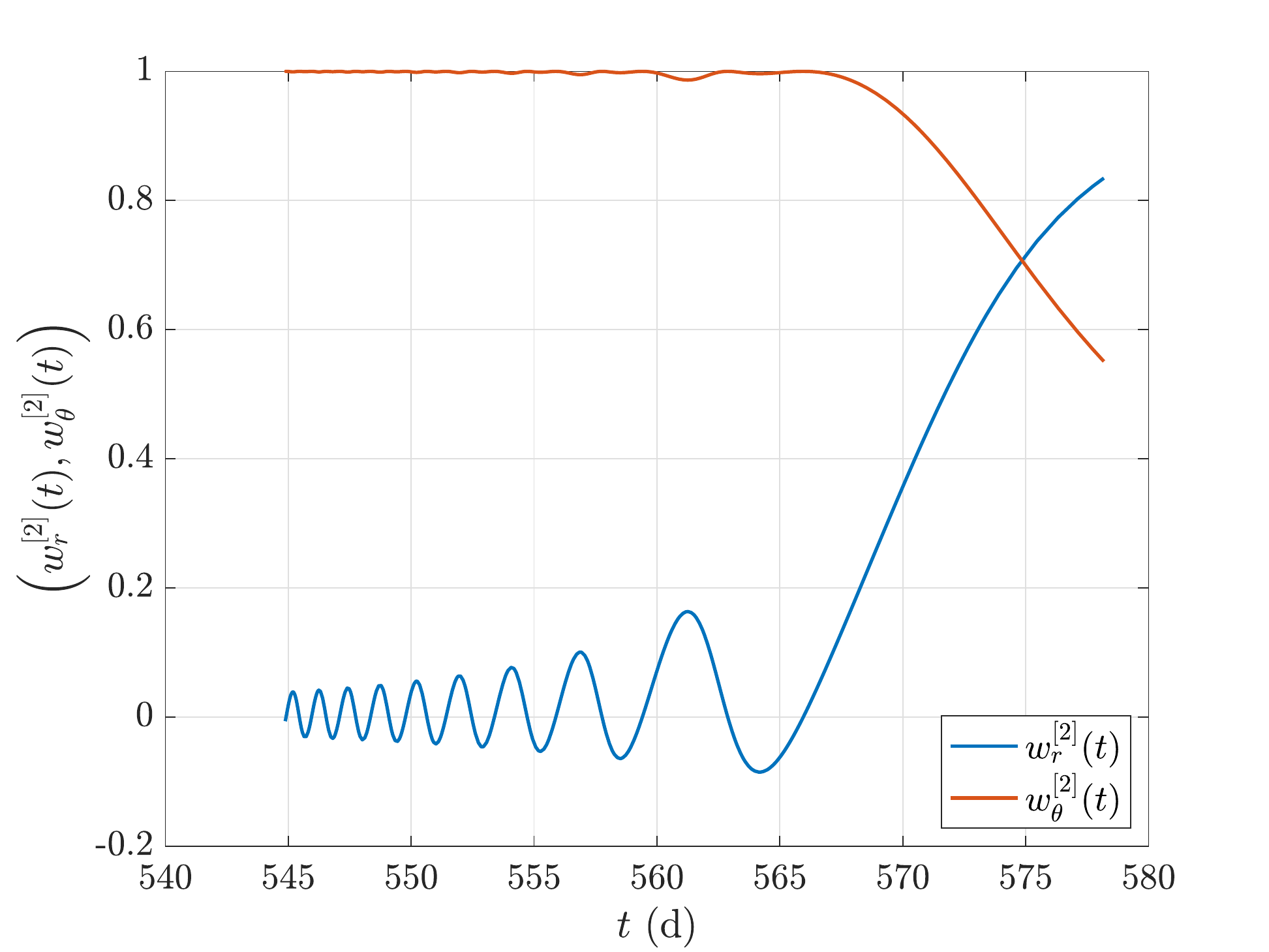}}
  \caption{Optimal trajectory and control for the Earth escape phase.\label{fig:pert-spacecraft_Phase_2}}
\end{figure}

Next, Fig~\ref{fig:pert-Trajectory_Phase_2} shows the optimized two-dimensional trajectory in Cartesian coordinates $\left(x^{[2]}(t),y^{[2]}(t)\right) = \left(r^{[2]}(t)\cos\left(\theta^{[2]}(t)\right),r^{[2]}(t)\sin\left(\theta^{[2]}(t)\right)\right)$ during phase 2 (Earth escape).  It is seen for phase 2 that the spacecraft starts in a circular orbit of radius $6.6 R_E$ relative to Earth and terminates a distance of $1.01R_E^{SOI}$ from the Earth.  Next, Fig.~\ref{fig:pert-eSC_Phase_2} shows the eccentricity of the spacecraft as a function of time.  It is seen that the spacecraft makes multiple revolutions around the Earth in phase 2 such that the eccentricity remains between zero and 0.2 for approximately 2/3 of the phase. Then, during the last 1/3 of the phase, the eccentricity quickly grows to larger than unity at which point the spacecraft orbit transitions from elliptic to hyperbolic relative to the Earth.  Escape from Earth (that is, an eccentricity that exceeds unity) occurs at approximately 75~percent of the way into the phase.  Finally, the eccentricity at the end of the phase is approximately two.  Next, Fig.~\ref{fig:pert-Control_Components_Phase_2} shows the radial and transverse components of the thrust direction, $w_r^{[2]}(t)$ and $w_\theta^{[2]}(t)$, respectively. During the portion of the phase where the eccentricity is small, the radial component of thrust remains close to zero while the transverse component of the thrust remains near unity.  It is noted that during this early part of the phase the radial component of the thrust oscillates with increasing amplitude about zero, thereby indicating that the spacecraft is being propelled further from the Earth.  Although the thrust points primarily in the tangential direction during the elliptic portion of the Earth escape phase, a small fraction of the thrust direction still points in the radial direction (see the nonzero $w_{r}^{[2]}$-component in Fig.~\ref{fig:pert-Control_Components_Phase_2}).  This small but nonzero radial component during the elliptic portion of the Earth escape phase keeps the spacecraft under the gravitational influence of the Earth for a longer duration while maintaining a larger velocity.  As a result, the transfer becomes more energy-efficient (thereby reducing the time required to complete the entire transfer) than it would be if the radial component of thrust during this segment was zero.  Then, during the hyperbolic segment of the phase, the radial component of thrust increases to approximately 0.9 while the transverse component of thrust decreases steadily to approximately 0.55.  The optimized duration of the Earth escape phase is $33.31~\textrm{d}$.

\begin{figure}[h]
  \centering
  \subfloat[Optimal two-dimensional trajectory.\label{fig:pert-Trajectory_Phase_3}]{\includegraphics[width=3.25in]{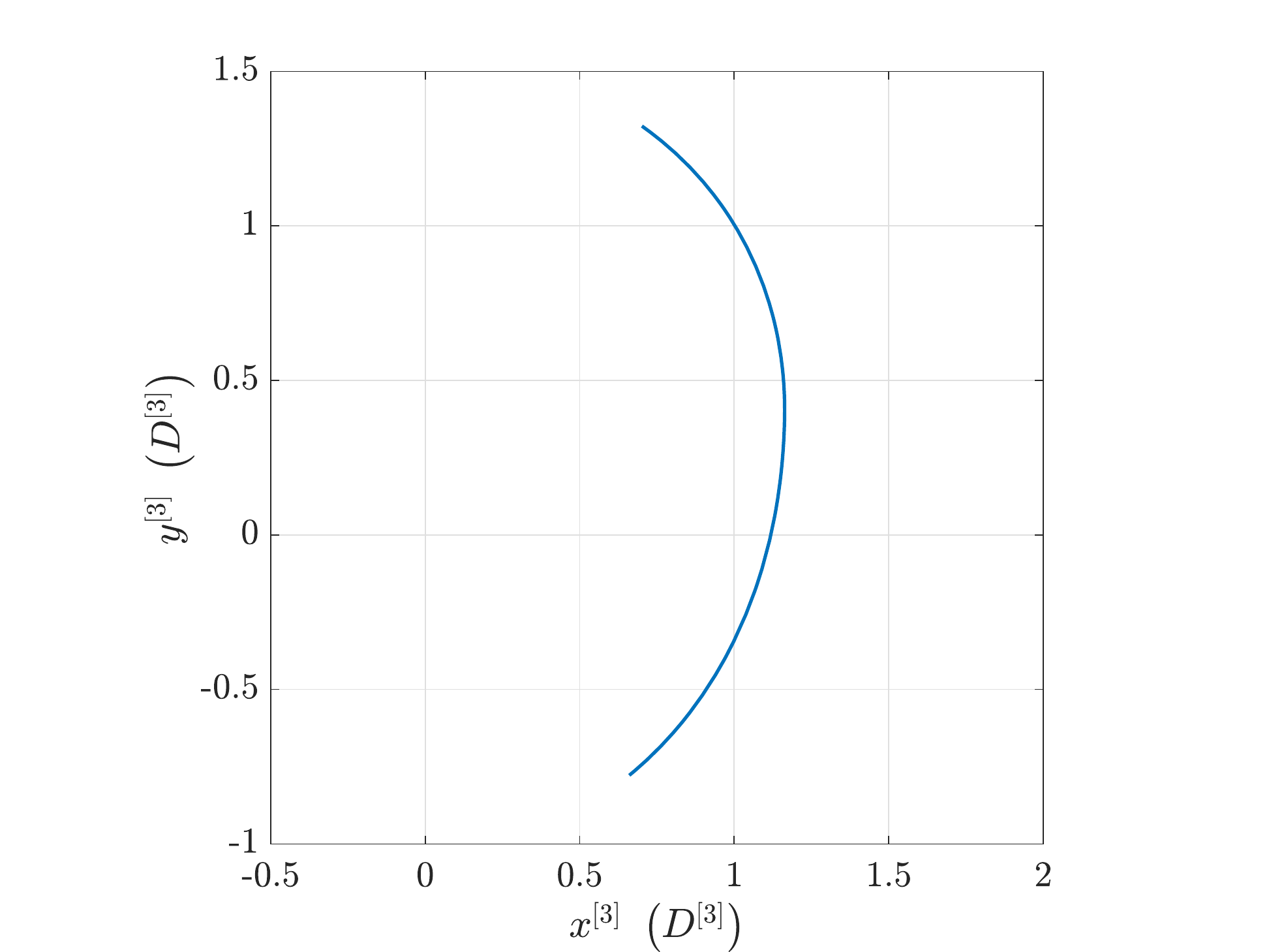}}~~\subfloat[Spacecraft eccentricity. \label{fig:pert-eSC_Phase_3}]{\includegraphics[width=3.25in]{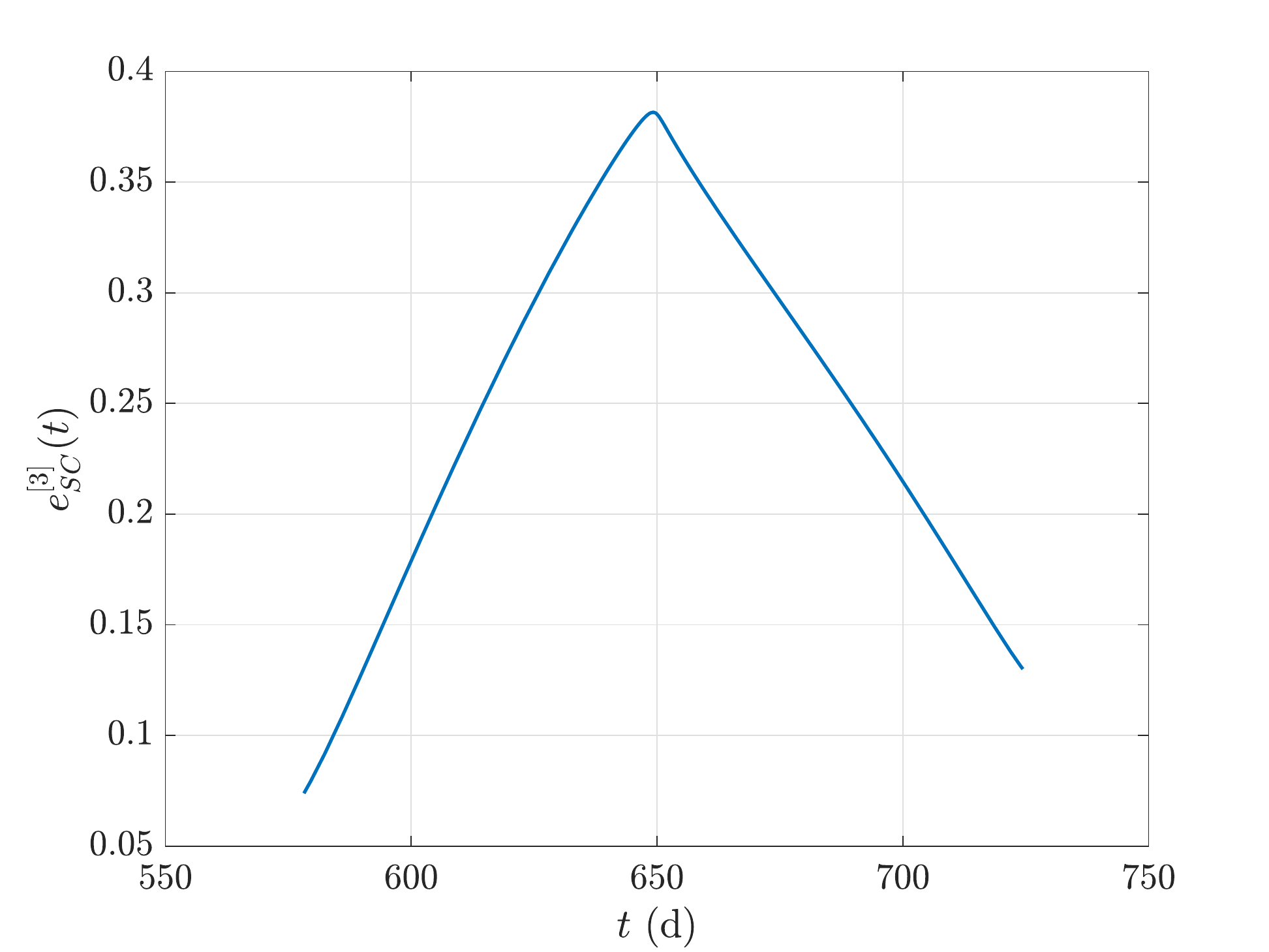}}

  \subfloat[Optimal controls $\left(w_r^{[3]},w_\theta^{[3]}\right)$ vs. $t$. \label{fig:pert-Control_Components_Phase_3}]{\includegraphics[width=3.25in]{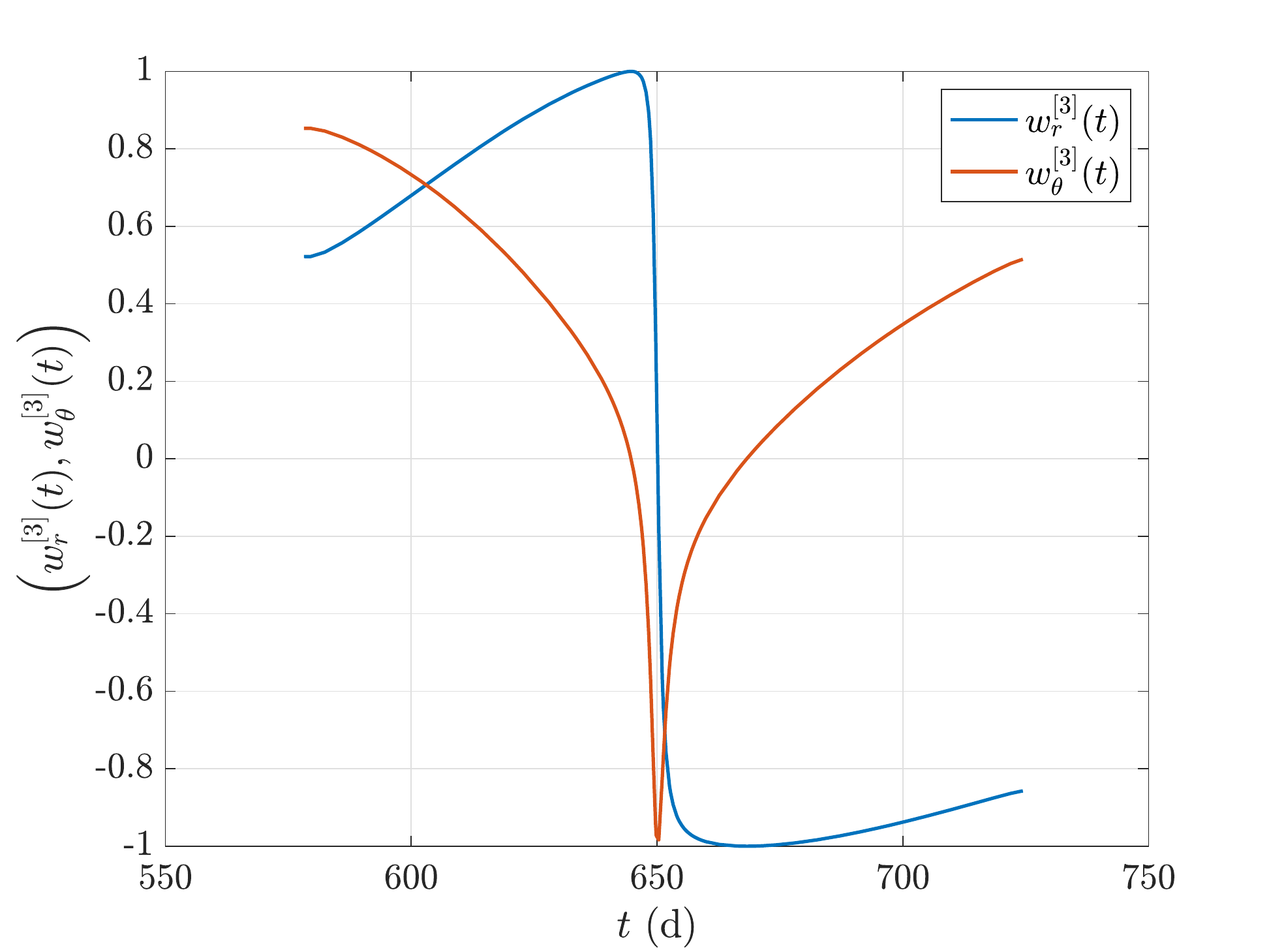}}
  \caption{Optimal trajectory and control for the heliocentric transfer phase.\label{fig:pert-spacecraft_Phase_3}}
\end{figure}

Next, Fig,.~\ref{fig:pert-Trajectory_Phase_3} shows the optimal two-dimensional trajectory in Cartesian coordinates $\left(x^{[3]}(t),y^{[3]}(t)\right) = \left(r^{[3]}(t)\cos\left(\theta^{[3]}(t)\right),r^{[3]}(t)\sin\left(\theta^{[3]}(t)\right)\right)$ during phase 3 (heliocentric transfer from Earth to Mars). It is seen that the spacecraft only makes a partial orbit around the Sun. Figure~\ref{fig:pert-eSC_Phase_3} shows the eccentricity of the spacecraft over the duration of phase 3 where it stays in a range of 0.07 to 0.37. Figure~\ref{fig:pert-Control_Components_Phase_3} shows the radial and transverse components of the thrust direction, $w_r^{[3]}(t)$ and $w_\theta^{[3]}(t)$, respectively.  It is seen that the thrust starts in a direction that enables the spacecraft to escape the gravitational field of the Earth.  Note, however, that approximately halfway through the interplanetary phase the thrust changes direction resulting in a retrograde maneuver.  The reason that the thrust direction eventually lies opposite the direction of motion is because the speed of the spacecraft must decrease in order to arrive in Mars orbit (because Mars is moving slower than Earth).   The duration of this phase is $146.25~\textrm{d}$.

\begin{figure}[h]
  \centering
  \subfloat[Optimal two-dimensional trajectory.\label{fig:pert-Trajectory_Phase_4}]{\includegraphics[width=3.25in]{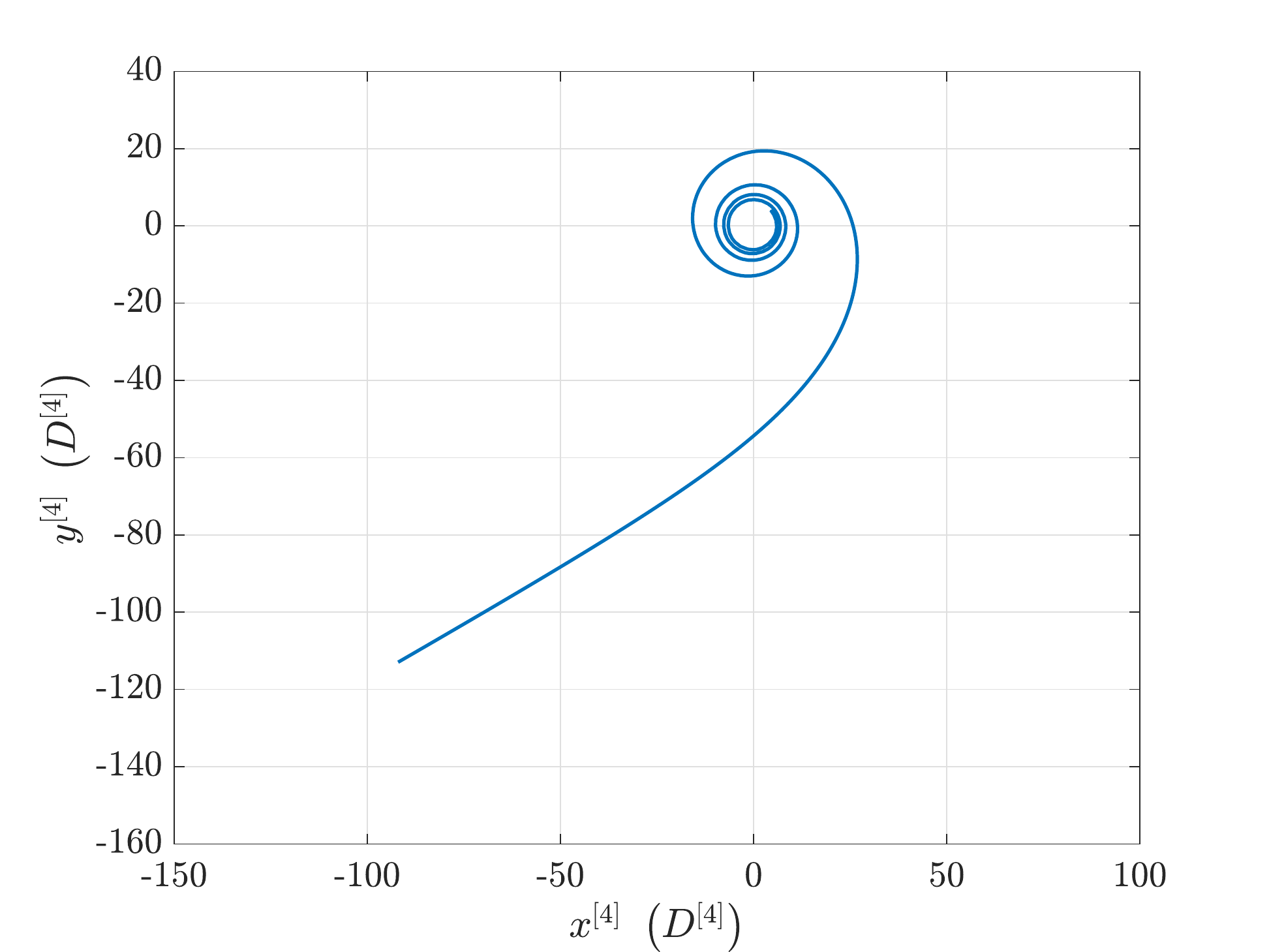}}~~\subfloat[Eccentricity of the spacecraft. \label{fig:pert-eSC_Phase_4}]{\includegraphics[width=3.25in]{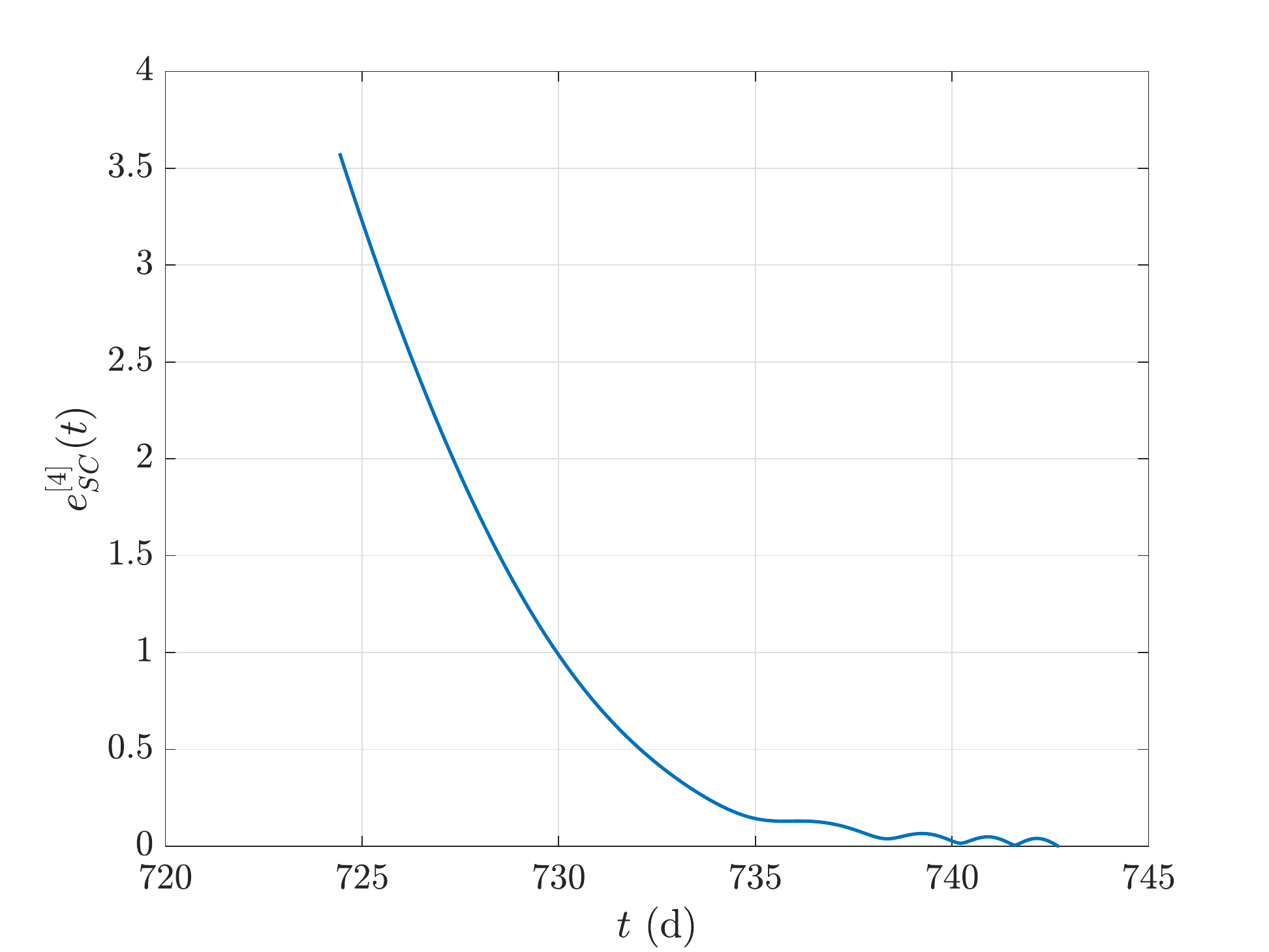}}

  \subfloat[Optimal controls $\left(w_r^{[4]},w_\theta^{[4]}\right)$ vs. $t$.\label{fig:pert-Control_Components_Phase_4}]{\includegraphics[width=3.25in]{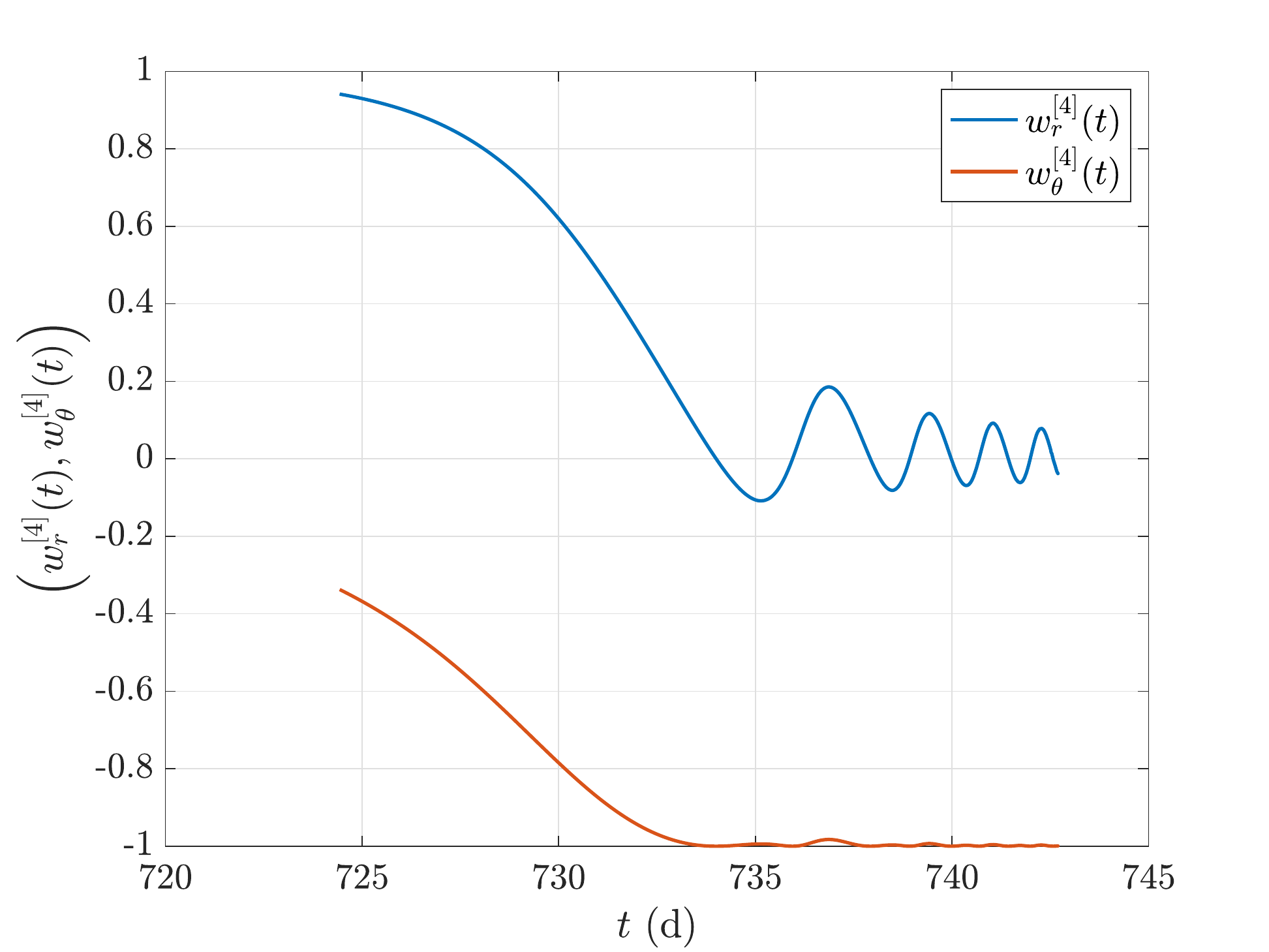}}

  \caption{Optimal trajectory and control for the Mars capture phase. \label{fig:pert-spacecraft_Phase_4}}
\end{figure}

Next, Fig.~\ref{fig:pert-Trajectory_Phase_4} shows the optimized two-dimensional trajectory of the spacecraft in Cartesian coordinates $\left(x^{[4]}(t),y^{[4]}(t)\right) = \left(r^{[4]}(t)\cos\left(\theta^{[4]}(t)\right),r^{[4]}(t)\sin\left(\theta^{[4]}(t)\right)\right)$ during phase 4 (Mars capture). It is seen for phase 4 that the spacecraft starts at a distance of $0.86R_M^{SOI}$ from Mars and terminates in a circular orbit of radius $6.0 R_M$ relative to Mars. Next, Fig.~\ref{fig:pert-eSC_Phase_4} shows the eccentricity of the spacecraft as a function of time. It is seen that the spacecraft starts hyperbolic relative to Mars for the first $1/4$ of the phase such that the eccentricity decreases from 3.5 to unity, where the spacecraft orbit transitions from hyperbolic to elliptic relative to Mars. Then, during the last $3/4$ of the phase, the spacecraft makes multiple revolutions around Mars such that the eccentricity decreases to 0, which corresponds to the terminal circular orbit. Capture to Mars (that is, an eccentricity that does not exceed unity) occurs at approximately 32~percent of the way into the phase. Next, Fig.~\ref{fig:pert-Control_Components_Phase_4} shows the radial and transverse components of the thrust direction, $w_r^{[4]}(t)$ and $w_\theta^{[4]}(t)$, respectively. During the hyperbolic segment of the phase, the radial component of thrust decreases from approximately 0.95 to 0.6 while the transverse component of thrust decreases from approximately -0.35 to -0.8. Then, during the portion of the phase where the eccentricity is small, the radial component of thrust remains close to zero while the transverse component of the thrust remains near -1. It is noted that during this later part of the phase the radial component of the thrust oscillates with decreasing amplitude about zero, thereby indicating that the spacecraft is being propelled towards Mars. The optimized duration of the Mars capture phase is $18.27~\textrm{d}$.

\begin{figure}[h!]
  \centering
  \includegraphics[width=3.25in]{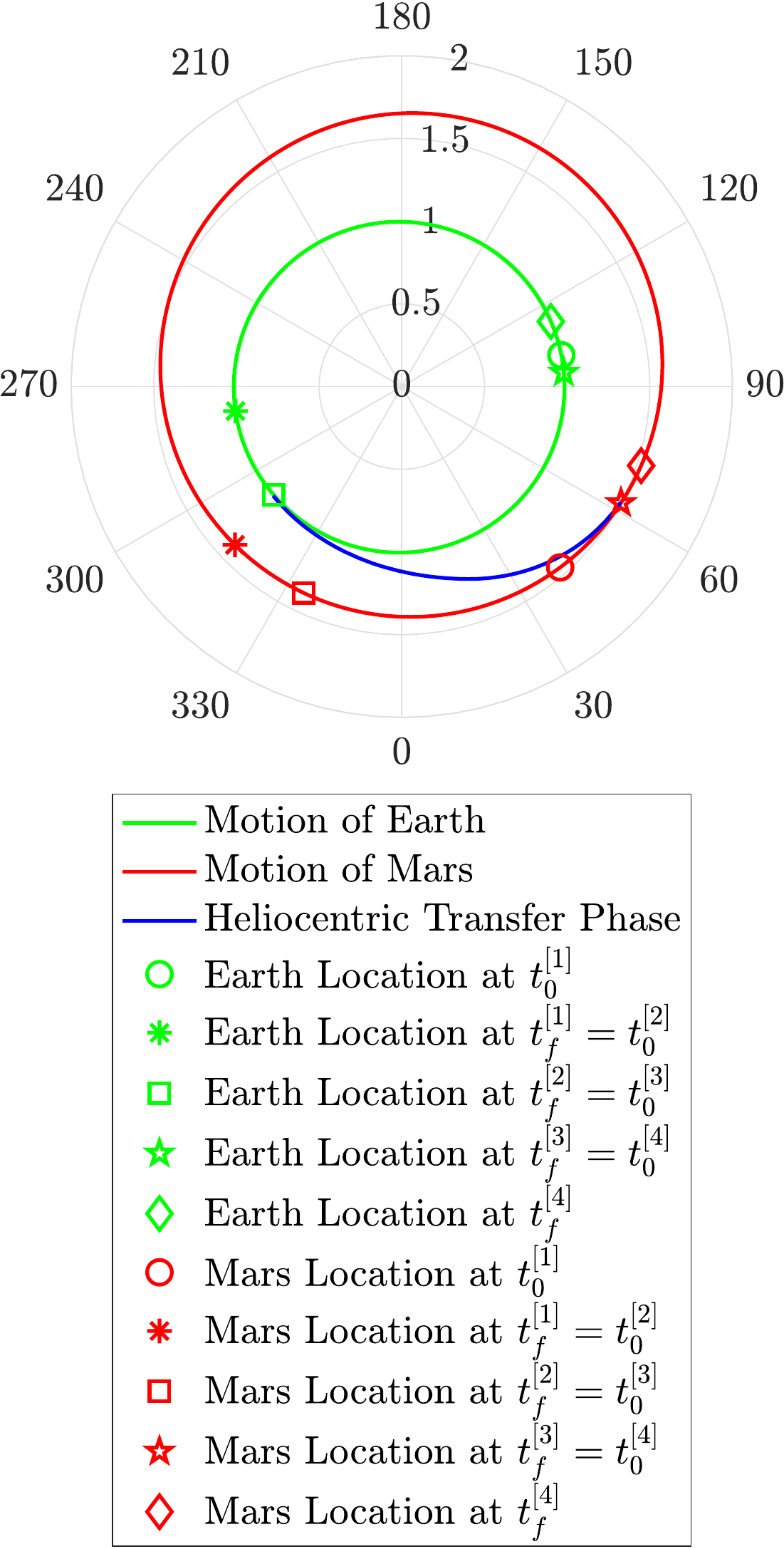}
  \caption{Four-phase minimum-time interplanetary orbit transfer.\label{fig:pert-Transfer_Solution}}
\end{figure}

Lastly, a higher-level view of the orbital transfer between Earth and Mars is shown in Fig.~\ref{fig:pert-Transfer_Solution}.  In particular, Fig.~\ref{fig:pert-Transfer_Solution} shows the motion of the planets and the spacecraft during the Earth escape, heliocentric, and Mars capture phases.  It is noted that, using a date of 1 January 2019 for the start of the planetary alignment phase, the interplanetary orbital transfer starts on 28 June 2020. The interplanetary transfer has a minimum-time solution of $197.83~\textrm{d}$. 

\clearpage

\section{Conclusions\label{sect:conclusions}}

A numerical optimization study of minimum-time low-thrust Earth-to-Mars orbital trajectory optimization was performed.  Two different cases of planetary motion were considered: circular and elliptic.  In addition, the problem was solved without and with third-body gravity perturbations. The trajectory optimization problem was formulated as a four-phase optimal control problem where the phases consisted of planetary alignment, Earth escape, heliocentric transfer, and Mars capture.  Because the problem was divided into four phases, interior-point event constraints were imposed to transform the variables at the terminus of one phase to the variables used at the start of the subsequent phase. It was found that the duration of the planetary alignment phases increased as a function of thrust specific force while the overall transfer time decreased.  For all values of thrust specific force considered, it was found that on average the duration of phases 1, 2, 3, and 4 for the elliptic with gravity perturbations case were, respectively, $545.99~\textrm{d}$, $33.53~\textrm{d}$, $143.42~\textrm{d}$, and $18.69~\textrm{d}$, while the average transfer time was $195.63~\textrm{d}$. Moreover, the departure dates from Earth for these values of thrust specific force were found to lie between 28 June 2020 and 1 July 2020.  In addition, the results for circular planetary motion were found to be qualitatively similar to those for elliptic planetary motion. Finally, it was found that, using the formulation developed in this study, the transfer times for the case of circular planetary motion without a planetary alignment phase (that is, a three-phase formulation) were lower than the transfer times obtained in a previous study.  

\section*{Acknowledgments}

The authors gratefully acknowledge support for this research from the U.S.~National Science Foundation under grant DMS-1819002 and from the National Aeronautics and Space Administration under grant NNX15AI10H through the University of Central Florida NASA Space Grant Consortium and Space Florida.  The authors would also like to thank one of the anonymous reviewers of this paper for identifying the reason for the existence of a negative radial component of thrust during the Earth escape phase.  

\clearpage

\renewcommand{\baselinestretch}{1}\normalsize\normalfont
% \bibliography{references.bib}

\appendix

\renewcommand{\baselinestretch}{2}\normalsize\normalfont

\section*{Appendix}

\section{Derivation of Coordinate Transformations\label{sect:coordinate-transformations}}

The derivations are now provided for the coordinate transformations from planet-centered coordinates to heliocentric coordinates.  Note that the results of this section can be applied to transform from either the Earth or the Mars coordinates to the heliocentric coordinates.  As a result, the derivation is provided only once for both coordinate transformations, leading to either Eq.~\eqref{events-earth-to-heliocentric} or \eqref{events-heliocentric-to-mars}.

\subsection{Notation, Conventions, and Coordinate Systems\label{sect:notation}}

\begin{figure}[h!]
  \centering
  \includegraphics[width=3.25in]{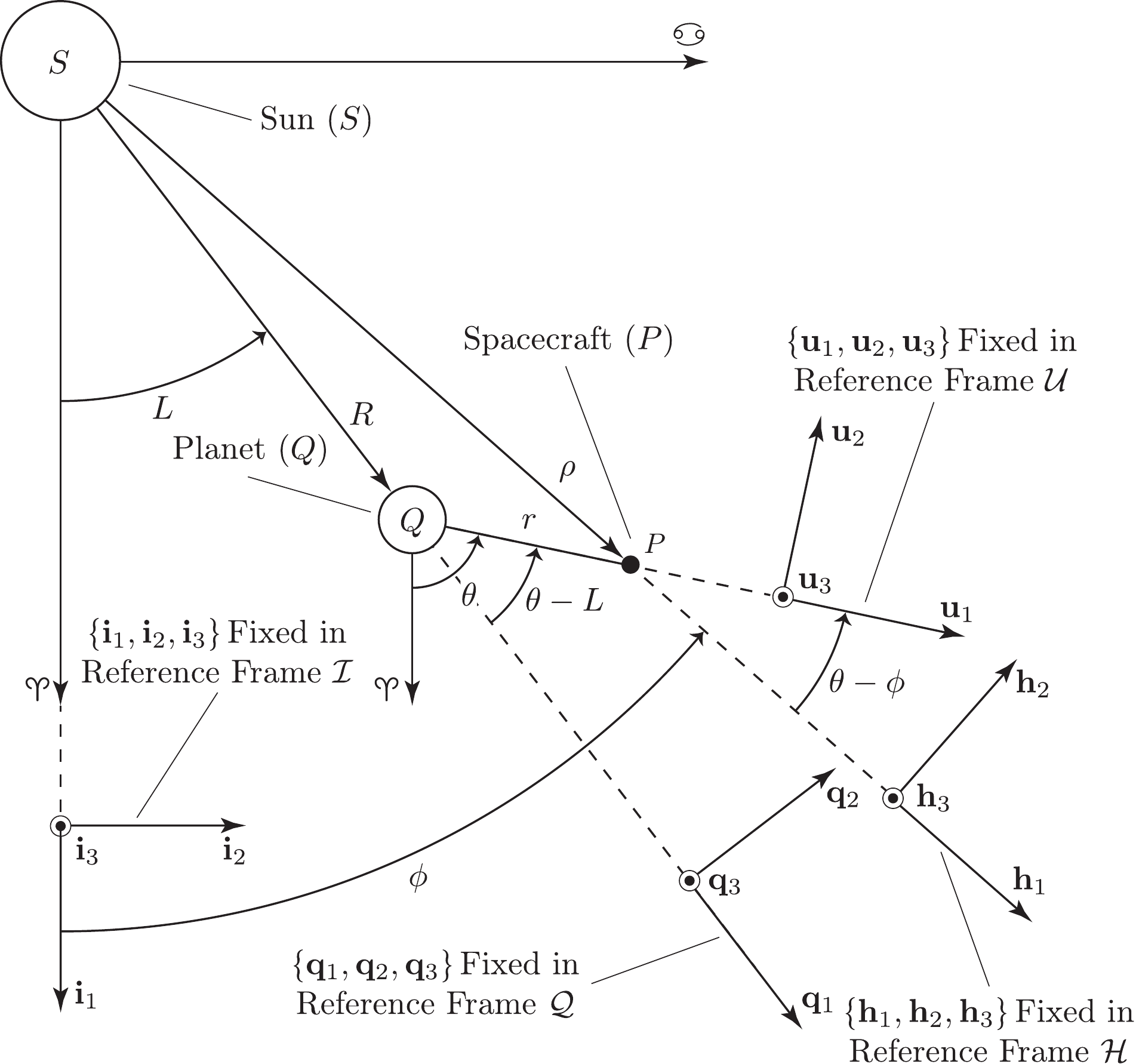}
  \caption{Schematic of Sun $(S)$, planet $(Q)$, spacecraft $(P)$ along with the variables $\rho$, $\phi$, $r$, and $\theta$ that define the position of the spacecraft relative to the Sun.\label{fig:sun-planet-spacecraft-generic}}
\end{figure}

Let $S$, $Q$, and $P$ denote, respectively, the locations of the Sun, the planet (central body), and the spacecraft.  Furthermore, assume that the Sun $(S)$ is an inertially fixed point.  Next, assume that $S$, $Q$, and $P$ all lie in the ecliptic plane where the ecliptic plane is assumed to be an inertial reference frame.  Also, assume that the $\m{N}=\m{U}_{\Aries}\times\m{U}_{\cancer}$ is the normal to the ecliptic plane, where $\m{U}_{\Aries}$ and $\m{U}_{\cancer}$ are the unit vectors along the First Point of Aries and the First Point of Cancer, respectively.  Finally, let $\m{U}_{SQ}$, $\m{U}_{SP}$, and $\m{U}_{QP}$, denote, respectively, the unit vectors in the directions from $S$ to $Q$, $S$ to $P$, and $Q$ to $P$.  The aforementioned unit vectors are then used to define the following reference frames:
\begin{enumerate}[1.]
\item $\C{I}=\left\{\m{U}_{\Aries},\m{N}\right\}$--plane = inertial reference frame.
\item $\C{Q}=\left\{\m{U}_{SQ},\m{N}\right\}$--plane = (Sun, central body) reference frame.
\item $\C{H}=\left\{\m{U}_{SP},\m{N}\right\}$--plane = (Sun, spacecraft) reference frame.
\item $\C{U}=\left\{\m{U}_{QP},\m{N}\right\}$--plane = (central body, spacecraft) reference frame.
\end{enumerate}
The following coordinate systems are then fixed, respectively, in reference frames $\C{I}$, $\C{Q}$, $\C{H}$, and $\C{U}$:
\begin{enumerate}[1.]
\item $\C{I}$-fixed coordinate system:  $(S, \m{i}_1, \m{i}_2, \m{i}_3)$
  \begin{equation}\label{inertially-fixed-coordinates}
    \begin{array}{lcl}
      \textrm{Origin:} & S, \\
      \m{i}_1 & = & \textrm{Along }\aries, \\
      \m{i}_3 & = & \textrm{Normal to Ecliptic Plane}, \\
      \m{i}_2 & = & \m{i}_3 \times \m{i}_1 . 
    \end{array}
  \end{equation} 
\item $\C{Q}$-fixed coordinate system:   $(S, \m{q}_1, \m{q}_2, \m{q}_3)$
  \begin{equation}\label{sun-earth-fixed-coordinates}
    \begin{array}{lcl}
      \textrm{Origin:} & S, \\
      \m{q}_1 & = & \textrm{Along }\m{U}_{SQ}, \\
      \m{q}_3 & = & \m{N}, \\
      \m{q}_2 & = & \m{q}_3 \times \m{q}_1.
    \end{array}
  \end{equation}
\item $\C{H}$-fixed coordinate system:  $(S, \m{h}_1, \m{h}_2, \m{h}_3)$
  \begin{equation}\label{sun-spacecraft-fixed-coordinates}
    \begin{array}{lcl}
      \textrm{Origin:} & S, \\
      \m{h}_1 & = & \textrm{Along }\m{U}_{SP}, \\
      \m{h}_3 & = & \m{N}, \\
      \m{h}_2 & = & \m{h}_3 \times \m{h}_1.  
    \end{array}
  \end{equation}
\item $\C{U}$-fixed coordinate system:    $(Q, \m{u}_1, \m{u}_2, \m{u}_3)$
  \begin{equation}\label{earth-spacecraft-fixed-coordinates}
    \begin{array}{lcl}
      \textrm{Origin:} & Q, \\
      \m{u}_1 & = & \textrm{Along }\m{U}_{QP}, \\
      \m{u}_3 & = & \m{N}, \\
      \m{u}_2 & = & \m{u}_3 \times \m{u}_1.
    \end{array}
  \end{equation}
\end{enumerate}
Then, using Fig.~\ref{fig:sun-planet-spacecraft-generic}, the bases $\left\{\m{q}_1,\m{q}_2,\m{q}_3\right\}$ and $\left\{\m{h}_1,\m{h}_2,\m{h}_3\right\}$  are related to the basis $\left\{\m{u}_1,\m{u}_2,\m{u}_3\right\}$, respectively, as
\begin{equation}\label{e2u}
  \begin{array}{lcl}
    \m{q}_1 & = & \cos\left(\theta-L\right)\m{u}_1 - \sin\left(\theta-L\right)\m{u}_2, \\
    \m{q}_2 & = & \sin\left(\theta-L\right)\m{u}_1 + \cos\left(\theta-L\right)\m{u}_2, \\
    \m{q}_3 & = & \m{u}_3,
  \end{array}
\end{equation}
and
\begin{equation}\label{h2u}
  \begin{array}{lcl}
    \m{h}_1 & = & \cos\left(\theta-\phi\right)\m{u}_1 - \sin\left(\theta-\phi\right)\m{u}_2, \\
    \m{h}_2 & = & \sin\left(\theta-\phi\right)\m{u}_1 + \cos\left(\theta-\phi\right)\m{u}_2, \\
    \m{h}_3 & = & \m{u}_3.
  \end{array}
\end{equation}
Using the aforementioned coordinate systems, the following variables are used to define the position and inertial velocity of the spacecraft relative to the Sun in terms of the coordinate systems $(S, \m{h}_1, \m{h}_2, \m{h}_3)$: 
\begin{equation}
  \begin{array}{lcl}
    \rho & = & \textrm{Distance from }S\textrm{ to }P, \\
    \phi & = & \textrm{Angle from }\m{i}_1\textrm{ to }\m{h}_1, \\
    v_\rho & = & \dot{\rho}, \\
    v_\phi & = & \rho\dot{\phi}.
  \end{array}
\end{equation}
Also, the following variables are used to define the position and inertial velocity of the spacecraft relative to the Sun in terms of the coordinate systems $(S, \m{q}_1, \m{q}_2, \m{q}_3)$ and $(Q, \m{u}_1, \m{u}_2, \m{u}_3)$:
\begin{equation}
  \begin{array}{lcl}
    r & = & \textrm{Distance from }Q\textrm{ to }P, \\
    \theta & = & \textrm{Angle from }\m{i}_1\textrm{ to }\m{u}_1, \\
    v_r & = & \dot{r}, \\
    v_\theta & = & r\dot{\theta}.
  \end{array}
\end{equation}
Finally, because the position of the planet is known, the quantities $R$ and $L$, and their rate of changes, $\dot{R}$ and $\dot{L}$, are considered to be known.  The variables $(\rho,\phi)$ and $(r,\theta)$ are then used to derive independent expressions for the position of the spacecraft relative to the Sun and the inertial velocity of the spacecraft.  Section \ref{sect:relationship-planet-centered-sun-centered} derives these independent expressions and arrives at relationships between $(\rho,\phi,v_\rho,v_\phi)$ and $(r,\theta,v_r,v_\theta)$.  A schematic of the various reference frames, coordinate systems, and variables is shown in Fig.~\ref{fig:sun-planet-spacecraft-generic}.

\subsection{Relationship Between Planet-Centered and Sun-Centered Coordinates\label{sect:relationship-planet-centered-sun-centered}}

The transformation from planet-centered to Sun-centered coordinates is now derived using Fig.~\ref{fig:sun-planet-spacecraft-generic}.  First, one expression for the position of the spacecraft relative to the Sun is expressed in the basis $\left\{\m{h}_1,\m{h}_2, \m{h}_3\right\}$ as 
\begin{equation}\label{position-spacecraft-h}
  \m{r}_{P/S} = \rho\m{h}_1.
\end{equation}
Now, because the angular velocity of reference frame $\m{H}$ as viewed by an observer in reference frame $\C{I}$ is $\lrsup{\C{I}}{\g{\omega}}{\C{H}}=\dot{\phi}\m{h}_3$, the inertial velocity of the spacecraft is obtained from the transport theorem \cite{Rao2006} as
\begin{equation}\label{inertial-velocity-spacecraft-h}
%   \begin{split}
    \lsup{\C{I}}{\m{v}_P}  = \lsupddt{\C{I}\hspace{-4pt}}{}\left(\m{r}_{P/S}\right) = \lsupddt{\C{H}\hspace{-4pt}}{}\left(\m{r}_{P/S}\right)  + \lrsup{\C{I}}{\g{\omega}}{\C{H}}\times \m{r}_{P/S}  = \dot{\rho}\m{h}_1 + \dot{\phi}\m{h}_3 \times \rho\m{h}_1  = \dot{\rho}\m{h}_1 + \rho\dot{\phi}\m{h}_2 \equiv  v_\rho\m{h}_1 + v_{\phi}\m{h}_2.
%   \end{split}
\end{equation}
Using expressions for $\m{h}_1$ and $\m{h}_2$ in terms of $\m{u}_1$ and $\m{u}_2$ given in Eq.~\eqref{h2u}, the inertial velocity of the spacecraft is expressed in the basis $\left\{\m{u}_1,\m{u}_2,\m{u}_3\right\}$ as
\begin{equation}
  \lsup{\C{I}}{\m{v}_P}  = v_\rho\left[\cos\left(\theta-\phi\right)\m{u}_1 - \sin\left(\theta-\phi\right)\m{u}_2\right] + v_{\phi}\left[\sin\left(\theta-\phi\right)\m{u}_1 + \cos\left(\theta-\phi\right)\m{u}_2\right].
\end{equation}
Next, a second expression for the position of the spacecraft relative to the Sun is given in terms of the bases $\left\{\m{q}_1,\m{q}_2,\m{q}_3\right\}$ and $\left\{\m{u}_1,\m{u}_2,\m{u}_3\right\}$ as
\begin{equation}\label{position-spacecraft-e-and-u}
  \m{r}_{P/S} = \m{r}_{Q/S} + \m{r}_{P/Q} = R\m{q}_1 + r\m{u}_1.
\end{equation}
The inertial velocity of the spacecraft can then be written as 
\begin{equation}
  \lsup{\C{I}}{\m{v}_P}  = \lsupddt{\C{I}\hspace{-4pt}}{}\left(\m{r}_{Q/S}\right) + \lsupddt{\C{I}\hspace{-4pt}}{}\left(\m{r}_{P/Q}\right) = \lsup{\C{I}}{\m{v}_Q} + \lsup{\C{I}}{\m{v}_{P/Q}}.
\end{equation}
Now, using the fact that the angular velocity of reference frame $\m{Q}$ as viewed by an observer in reference frame $\C{I}$ is $\lrsup{\C{I}}{\g{\omega}}{\C{Q}}=\dot{L}\m{q}_3$, the inertial velocity of point $Q$ is obtained from the transport theorem \cite{Rao2006} as
\begin{equation}\label{inertial-velocity-earth-e}
  \lsup{\C{I}}{\m{v}_Q} = \lsupddt{\C{I}\hspace{-4pt}}{}\left(\m{r}_{Q/S}\right) = \lsupddt{\C{Q}\hspace{-4pt}}{}\left(\m{r}_{Q/S}\right) + \lrsup{\C{I}}{\g{\omega}}{\C{Q}}\times \m{r}_{Q/S} = \dot{R}\m{q}_1 + \dot{L}\m{q}_3 \times  R\m{q}_1 = \dot{R}\m{q}_1 + R\dot{L}\m{q}_2.
\end{equation}
Furthermore, because the angular velocity of reference frame $\m{U}$ as viewed by an observer in reference frame $\C{I}$ is $\lrsup{\C{I}}{\g{\omega}}{\C{U}}=\dot{\theta}\m{u}_3$, the inertial velocity of the spacecraft relative to the central body is obtained from the transport theorem \cite{Rao2006} as
\begin{equation}\label{inertial-velocity-spacecraft-relative-to-earth-u}
  \lsup{\C{I}}{\m{v}_{P/Q}} = \lsupddt{\C{I}\hspace{-4pt}}{}\left(\m{r}_{P/Q}\right) = \lsupddt{\C{U}\hspace{-4pt}}{}\left(\m{r}_{P/Q}\right) + \lrsup{\C{I}}{\g{\omega}}{\C{U}}\times \m{r}_{P/Q} = \dot{r}\m{u}_1 + \dot{\theta}\m{u}_3 \times  r\m{u}_1 \equiv \dot{r}\m{u}_1 + r\dot{\theta}\m{u}_2 \equiv v_r\m{u}_1 + v_\theta\m{u}_2.
\end{equation}
Adding the results of Eqs.~\eqref{inertial-velocity-earth-e} and \eqref{inertial-velocity-spacecraft-relative-to-earth-u}, the second expression for the inertial velocity of the spacecraft is given as
\begin{equation}\label{inertial-velocity-spacecraft-e-and-u}
  \lsup{\C{I}}{\m{v}_P} = \lsup{\C{I}}{\m{v}_Q} + \lsup{\C{I}}{\m{v}_{P/Q}} = \dot{R}\m{q}_1 + R\dot{L}\m{q}_2 + \dot{r}\m{u}_1 + r\dot{\theta}\m{u}_2 = \dot{R}\m{q}_1 + R\dot{L}\m{q}_2 + v_r\m{u}_1 + v_{\theta}\m{u}_2.
\end{equation}

It is seen that Eqs.~\eqref{position-spacecraft-h} and \eqref{position-spacecraft-e-and-u} are two expressions for the position of the spacecraft relative to the Sun, while Eqs.~\eqref{inertial-velocity-spacecraft-h} and \eqref{inertial-velocity-spacecraft-e-and-u} are two expressions for the inertial velocity of the spacecraft.  Suppose now that the expressions given in Eqs.~\eqref{position-spacecraft-h} and \eqref{position-spacecraft-e-and-u} are projected into the basis $\left\{\m{u}_1,\m{u}_2,\m{u}_3\right\}$.  First, projecting the expression given in Eq.~\eqref{position-spacecraft-h} into the basis $\left\{\m{u}_1,\m{u}_2,\m{u}_3\right\}$ gives
\begin{equation}\label{projection-spacecraft-position-h-into-u}
  \begin{array}{lcl}
    \m{r}_{P/S}\cdot\m{u}_1 & = & \phantom{-}\rho\cos\left(\theta-\phi\right), \\
    \m{r}_{P/S}\cdot\m{u}_2 & = & -\rho\sin\left(\theta-\phi\right).
  \end{array}
\end{equation} 
Next, projecting the expression given in Eq.~\eqref{position-spacecraft-e-and-u} into the basis $\left\{\m{u}_1,\m{u}_2,\m{u}_3\right\}$ gives
\begin{equation}\label{projection-spacecraft-position-e-and-u-into-u}
  \begin{array}{lcl}
    \m{r}_{P/S}\cdot\m{u}_1 & = & \phantom{-}R\cos\left(\theta-L\right) + r, \\
    \m{r}_{P/S}\cdot\m{u}_2 & = & -R\sin\left(\theta-L\right). 
  \end{array}
\end{equation}
Setting the first and second expressions in Eq.~\eqref{projection-spacecraft-position-h-into-u} equal to the first and second expressions in Eq.~\eqref{projection-spacecraft-position-e-and-u-into-u}, respectively, and rearranging gives
\begin{equation}\label{coordinate-transformation-position}
  \begin{array}{lcl}
    r & = & -R\cos\left(\theta-L\right) + \rho\cos\left(\theta-\phi\right), \\
    0 & = &\phantom{-}R\sin(\theta-L) - \rho\sin\left(\theta-\phi\right).
  \end{array}
\end{equation}
Next, projecting the expression given in Eq.~\eqref{inertial-velocity-spacecraft-h} into the basis $\left\{\m{u}_1,\m{u}_2,\m{u}_3\right\}$ gives
\begin{equation}\label{inertial-velocity-h-projection-u}
  \begin{array}{lcl}
    \lsup{\C{I}}{\m{v}_P}\cdot\m{u}_1 & = & \phantom{-}v_\rho\cos(\theta-\phi) + v_\phi\sin(\theta-\phi), \\
    \lsup{\C{I}}{\m{v}_P}\cdot\m{u}_2 & = & -v_\rho\sin(\theta-\phi) + v_\phi\cos(\theta-\phi). \\
  \end{array}
\end{equation}
Then, projecting the expression given in Eq.~\eqref{inertial-velocity-spacecraft-e-and-u} into the basis $\left\{\m{u}_1,\m{u}_2,\m{u}_3\right\}$ gives
\begin{equation}\label{inertial-velocity-e-and-u-projection-u}
  \begin{array}{lcl}
    \lsup{\C{I}}{\m{v}_P}\cdot\m{u}_1 & = & \phantom{-}\dot{R}\cos(\theta-L) + R\dot{L}\sin(\theta-L) + v_r, \\
    \lsup{\C{I}}{\m{v}_P}\cdot\m{u}_2 & = & -\dot{R}\sin(\theta-L) + R\dot{L}\cos(\theta-L) + v_\theta.
  \end{array}
\end{equation}
Setting the first and second expressions in Eq.~\eqref{inertial-velocity-h-projection-u} equal to the first and second expressions in Eq.~\eqref{inertial-velocity-e-and-u-projection-u}, respectively, and rearranging gives
\begin{equation}\label{coordinate-transformation-velocity}
  \begin{array}{lcl}
    v_r & = & -\dot{R}\cos(\theta-L) - R\dot{L}\sin(\theta-L) + v_\rho\cos(\theta-\phi) + v_\phi\sin(\theta-\phi), \\
    v_\theta & = & \phantom{-}\dot{R}\sin(\theta-L) - R\dot{L}\cos(\theta-L) -v_\rho\sin(\theta-\phi) + v_\phi\cos(\theta-\phi).
  \end{array}
\end{equation}
Equation~\eqref{coordinate-transformation-position} provides the relationship between the variables $(\rho,\phi)$ and $(r,\theta)$, while Eq.~\eqref{coordinate-transformation-velocity} provides the relationship between the variables $(r,\theta,v_r,v_\theta)$ and the variables $(\rho,\phi,v_\rho,v_\phi)$.  Equations~\eqref{coordinate-transformation-position} and \eqref{coordinate-transformation-velocity} are used to develop the event constraints of Eqs.~\eqref{events-earth-to-heliocentric} and \eqref{events-heliocentric-to-mars}.  It is noted that Eq.~\eqref{events-earth-to-heliocentric} connects the terminus of the Earth escape phase to the start of the heliocentric phase, while Eq.~\eqref{events-heliocentric-to-mars} connects the terminus of the heliocentric phase to the start of the Mars capture phase.  When using Eqs.~\eqref{coordinate-transformation-position} and \eqref{coordinate-transformation-velocity} in the context of the event constraints given in Section \ref{sect:event-constraints}, the variables $(\rho,\phi,v_\rho,v_\phi)$ and $(r,\theta,v_r,v_\theta)$ have been replaced with the appropriate variables given in Table~\ref{tab:variables}.  

\section{Derivation of Gravity Perturbations\label{sect:perturbations}}

The derivations are now provided for the gravity perturbations given in Eq.~\eqref{gravity-perturbations}.  Note that for the Earth escape, Mars capture, and heliocentric phases (that is, phases 2, 3, and 4) Cowell's method~\cite{Bate1971} is used.  On the other hand, no gravity perturbations are included in phase 1 because the motion of the spacecraft is not included until phase 2. 

\subsection{Gravity Perturbations During Phases 2 and 4\label{sect:perturbations2,4}}

During phases 2 and 4, where the central bodies are Earth and Mars, respectively, the Sun is the only perturbing body.  Because the perturbing acceleration due to the Sun lies along $\m{h}_1$ (see Fig.~\ref{fig:sun-planet-spacecraft-generic}) while the dynamics are expressed in the basis $\left\{\m{u}_1, \m{u}_2, \m{u}_3\right\}$, the perturbing acceleration due to the Sun in the Earth escape and Mars captures phases is projected into the basis $\left\{\m{u}_1, \m{u}_2, \m{u}_3\right\}$.   % the unit vector $\m{h}_1$ needs to be obtained in terms of known quantities according to the phase.

Using Fig.~\ref{fig:sun-planet-spacecraft-generic}, the position of the spacecraft relative to the Sun is given as
\begin{equation}\label{r_P/S}
\m{r}_{P/S} = \m{r}_{Q/S} + \m{r}_{P/Q} = R\m{q}_1 + r\m{u}_1 = \rho\m{h}_1.
\end{equation}
% Therefore, 
% \begin{equation}\label{r_P/SRewritten}
%   \m{r}_{P/S} = \displaystyle \rho\m{h}_1.
% \end{equation}
Expressing $\m{q}_1$ in the basis $\{\m{u}_1, \m{u}_2, \m{u}_3\}$ gives
\begin{equation} \label{r_P/SConversion}
  \begin{split}
    \rho\m{h}_1 & = \displaystyle R\left[\cos\left(\theta-L_Q\right)\m{u}_1 - \sin\left(\theta-L_Q\right)\m{u}_2\right] + r\m{u}_1 \\
    & = \displaystyle \left[r + R\cos\left(\theta-L_Q\right)\right]\m{u}_1 - R\sin\left(\theta-L_Q\right)\m{u}_2,
  \end{split}
\end{equation}
where
\begin{equation}\label{h1-in-terms-of-u1-and-u2}
  \m{h}_1 = \displaystyle \frac{\left[r + R\cos\left(\theta-L_Q\right)\right]\m{u}_1 - R\sin\left(\theta-L_Q\right)\m{u}_2}{\rho} .
\end{equation}
Next, using the first expression in Eq.~\eqref{r_P/S}, $\m{h}_1$ can be written as
\begin{equation}\label{h_1UnitVectorRewritten}
    \m{h}_1 = \frac{\m{r}_{P/S}}{\|\m{r}_{P/S}\|} = \displaystyle \frac{\left[r + R\cos\left(\theta-L_Q\right)\right]\m{u}_1 - R\sin\left(\theta-L_Q\right)\m{u}_2}{\sqrt{\left(r + R\cos\left(\theta-L_Q\right)\right)^2 + \left(-R\sin\left(\theta-L_Q\right)\right)^2}}.
\end{equation}
Setting the results of Eqs.~\eqref{h1-in-terms-of-u1-and-u2} and~\eqref{h_1UnitVectorRewritten} equal to one another gives
\begin{equation}\label{VariableToSimplifyh_1}
\rho = \displaystyle \sqrt{\left(r + R\cos\left(\theta-L_Q\right)\right)^2 + \left(-R\sin\left(\theta-L_Q\right)\right)^2}  = \displaystyle \sqrt{r^2 + R^2 + 2rR\cos\left(\theta-L_Q\right)}.
\end{equation}
Then, using Cowell's method \cite{Bate1971}, the gravity perturbation of the Sun for use in phases 2 and 4 is given as
\begin{equation}\label{BasicCowell'sTerm}
\lsup{\C{I}}{\m{a}_p} = \displaystyle -\mu_S\left[\frac{\m{r}_{P/S}}{\|\m{r}_{P/S}\|^3} - \frac{\m{r}_{Q/S}}{\|\m{r}_{Q/S}\|^3}\right] =  -\mu_S\left[\frac{\rho}{\rho^3}\m{h}_1 - \frac{R}{R^3}\m{q}_1\right] = \displaystyle -\mu_S\left[\frac{1}{\rho^2}\m{h}_1 - \frac{1}{R^2}\m{q}_1\right].
\end{equation}
Now, using Eqs.~\eqref{h1-in-terms-of-u1-and-u2} and \eqref{e2u}, $\m{a}_p$ can be expressed in the basis $\left\{\m{u}_1,\m{u}_2,\m{u}_3\right\}$ as
\begin{equation}\label{Cowell'sTermInU}
\lsup{\C{I}}{\m{a}_p} = \displaystyle -\mu_S\left[\frac{\left[\left(r + R\cos\left(\theta - L_Q\right)\right)\m{u}_1 - R\sin\left(\theta - L_Q\right)\m{u}_2\right]}{\rho^3} - \frac{\left[\cos\left(\theta-L_Q\right)\m{u}_1 - \sin\left(\theta-L_Q\right)\m{u}_2\right]}{R^2}\right].
\end{equation}
Projecting $\lsup{\C{I}}{\m{a}_p}$ into the basis $\{\m{u}_1, \m{u}_2, \m{u}_3\}$ gives
\begin{equation}\label{Projectinga_P}
  \begin{array}{lcl}
    a_{pr} = \displaystyle \lsup{\C{I}}{\m{a}_p}\cdot\m{u}_1, \\
    a_{p\theta} = \displaystyle \lsup{\C{I}}{\m{a}_p}\cdot\m{u}_2.  
  \end{array}
\end{equation}
The quantities $a_{pr}$ and $a_{p\theta}$ are then given as
\begin{equation}\label{a_PTerms}
\begin{array}{lcl}
a_{pr} = \displaystyle -\mu_S\left[\frac{r}{\rho^3} + \left(\frac{R}{\rho^3} - \frac{1}{R^2}\right)\cos(\theta - L_Q)\right], \\
a_{p\theta} = \displaystyle -\mu_S\left(-\frac{R}{\rho^3} + \frac{1}{R^2}\right)\sin(\theta - L_Q). 
\end{array}
\end{equation}
It is noted again that the result given in Eq.~\eqref{a_PTerms} can be used in both the Earth escape phase (that is, phase 2) and the Mars capture phase (that is, phase 4).

\subsection{Gravity Perturbations During Phase 3\label{sect:perturbations3}}

During phase 3, where the central body is the Sun, both Earth and Mars are perturbing bodies. It is noted that the derivation shown below considers a generic perturbing body which could be either Earth or Mars.  The total perturbing acceleration for the heliocentric phase, which includes perturbations from both Earth and Mars, is given in \ref{sect:SimplifiedGravityPerturbations}. Because the perturbing acceleration due to a perturbing body $Q$ lies along $\m{u}_1$ (see Fig.~\ref{fig:sun-planet-spacecraft-generic}) while the dynamics are expressed in the basis $\{\m{h}_1, \m{h}_2, \m{h}_3\}$, the perturbing acceleration due to $Q$ in the heliocentric transfer phase is projected into the basis $\{\m{h}_1, \m{h}_2, \m{h}_3\}$.

Using Fig.~\ref{fig:sun-planet-spacecraft-generic}, the basis $\{\m{q}_1, \m{q}_2, \m{q}_3\}$ is related to the basis $\{\m{h}_1, \m{h}_2, \m{h}_3\}$ as
\begin{equation}\label{q2h}
\begin{array}{lcl}
\m{q}_1 & = & \cos\left(\phi-L_Q\right)\m{h}_1 - \sin\left(\phi-L_Q\right)\m{h}_2, \\

\m{q}_2 & = & \sin\left(\phi-L_Q\right)\m{h}_1 + \cos\left(\phi-L_Q\right)\m{h}_2, \\

\m{q}_3 & = & \m{h}_3.
\end{array}
\end{equation}
Next, from Eq.~\eqref{r_P/S}, the position of the spacecraft relative to the planet is given as
\begin{equation}\label{r_P/QRewritten}
\m{r}_{P/Q} = \displaystyle \m{r}_{P/S} - \m{r}_{Q/S} = \displaystyle \rho\m{h}_1 - R\m{q}_1 = r\m{u}_1.
\end{equation} 
Expressing $\m{q}_1$ in the basis $\{\m{h}_1, \m{h}_2, \m{h}_3\}$ using Eq.~\eqref{q2h} gives
\begin{equation}\label{r_P/QConversion}
\begin{split}
r\m{u}_1 & = \displaystyle \rho\m{h}_1 - R\left[\cos\left(\phi-L_Q\right)\m{h}_1 - \sin\left(\phi-L_Q\right)\m{h}_2\right] \\
& = \displaystyle \left[\rho - R\cos\left(\phi-L_Q\right)\right]\m{h}_1 + R\sin\left(\phi-L_Q\right)\m{h}_2,
\end{split}
\end{equation}
where
\begin{equation}\label{u1-in-terms-of-h1-and-h2}
\m{u}_1 = \displaystyle \frac{\left[\rho - R\cos\left(\phi-L_Q\right)\right]\m{h}_1 + R\sin\left(\phi-L_Q\right)\m{h}_2}{r}.
\end{equation}
Next, using the first expression in Eq.~\eqref{r_P/QRewritten}, $\m{u}_1$ can be written as
\begin{equation}\label{u_1UnitVectorRewritten}
\m{u}_1 = \displaystyle \frac{\m{r}_{P/Q}}{\|\m{r}_{P/Q}\|}  = \displaystyle \frac{\left[\rho - R\cos\left(\phi-L_Q\right)\right]\m{h}_1 
+ R\sin\left(\phi-L_Q\right)\m{h}_2}{\sqrt{\left(\rho - R\cos\left(\phi-L_Q\right)\right)^2 + \left(R\sin\left(\phi-L_Q\right)\right)^2}}.
\end{equation}
Setting the results of Eqs.~\eqref{u1-in-terms-of-h1-and-h2} and~\eqref{u_1UnitVectorRewritten} equal gives
\begin{equation}\label{VariableToSimplfyu_1}
\begin{array}{lcl}
r = \displaystyle \sqrt{\left(\rho - R\cos\left(\phi-L_Q\right)\right)^2 + \left(R\sin\left(\phi-L_Q\right)\right)^2} = \displaystyle \sqrt{\rho^2 + R^2 - 2R\rho\cos\left(\phi-L_Q\right)} .
\end{array}
\end{equation}
Then, using Cowell's method \cite{Bate1971}, the gravity perturbation due to $Q$ in phase 3 is given as 
\begin{equation}\label{BasicCowell'sTermHelio}
\lsup{\C{I}}{\m{a}_p} = \displaystyle -\mu_Q\left[\frac{\m{r}_{P/Q}}{\|\m{r}_{P/Q}\|^3} \right] = \displaystyle -\mu_Q\left[\frac{r}{r^3}\m{u}_1\right] = \displaystyle -\mu_Q\left[\frac{1}{r^2}\m{u}_1\right].
\end{equation}
Now, using Eqs.~\eqref{q2h} and \eqref{u1-in-terms-of-h1-and-h2}, $\lsup{\C{I}}{\m{a}_p}$ can be expressed in the basis $\left\{\m{h}_1,\m{h}_2,\m{h}_3\right\}$ as
\begin{equation}\label{Cowell'sTermInH}
\lsup{\C{I}}{\m{a}_p} = \displaystyle -\mu_Q\left[\frac{\left[\left(\rho - R\cos\left(\phi - L_Q\right)\right)\m{h}_1 + R\sin\left(\phi - L_Q\right)\m{h}_2\right]}{r^3}\right].
\end{equation}
Projecting $\lsup{\C{I}}{\m{a}_p}$ into the basis $\{\m{h}_1, \m{h}_2, \m{h}_3\}$ gives
\begin{equation}\label{Projectinga_PinH}
\begin{array}{lcl}
a_{pr} = \displaystyle \lsup{\C{I}}{\m{a}_p}\cdot\m{h}_1, \\
a_{p\theta} = \displaystyle \lsup{\C{I}}{\m{a}_p}\cdot\m{h}_2.
\end{array}
\end{equation}
The quantities $a_{pr}$ and $a_{p\theta}$ are then given as
\begin{equation}\label{a_PTermsHelio}
\begin{array}{lcl}
a_{pr} = \displaystyle -\mu_Q\left[\frac{\rho}{r^3} - \frac{R}{r^3}\cos\left(\phi - L_Q\right)\right], \\
a_{p\theta} = \displaystyle -\mu_Q\frac{R}{r^3}\sin(\phi-L_Q).
\end{array}
\end{equation}
It is noted again that the result given in Eq.~\eqref{a_PTermsHelio} can be used to obtain the perturbing acceleration of either Earth or Mars.

\subsection{Simplified Gravity Perturbations Using Cowell's Method\label{sect:SimplifiedGravityPerturbations}}

Using Cowell's method in Ref.~\cite{Bate1971}, the radial and transverse components of the gravity perturbations for phases 2, 3, and 4 are given as follows: 
\begin{equation}\label{SimplifiedGravityPerturbationsEq}
\begin{array}{lcl}
a_{pr}^{[2]} & = & \displaystyle -\mu_S^{[2]}\left[\frac{r^{[2]}}{r_{psE}^3} + \left(\frac{r_E^{[2]}}{r_{psE}^3}-\frac{1}{\left(r_E^{[2]}\right)^2}\right)\cos\left(\theta^{[2]}-L_E^{[2]}\right)\right], \vspace{5pt}\\

a_{p\theta}^{[2]} & = & \displaystyle - \mu_S^{[2]}\left(-\frac{r_E^{[2]}}{r_{psE}^3}+\frac{1}{\left(r_E^{[2]}\right)^2}\right)\sin\left(\theta^{[2]}-L_E^{[2]}\right), \vspace{5pt}\\

a_{pr}^{[3]} & = & \displaystyle -\mu_E^{[3]}\left[\frac{r^{[3]}}{r_{pe}^3} - \frac{r_E^{[3]}}{r_{pe}^3}\cos\left(\theta^{[3]}-L_E^{[3]}\right)\right] -\mu_M^{[3]}\left[\frac{r^{[3]}}{r_{pm}^3} - \frac{r_M^{[3]}}{r_{pm}^3}\cos\left(\theta^{[3]}-L_M^{[3]}\right)\right], \vspace{5pt}\\

a_{p\theta}^{[3]} & = & \displaystyle - \mu_E^{[3]}\frac{r_E^{[3]}}{r_{pe}^3}\sin\left(\theta^{[3]}-L_E^{[3]}\right) - \mu_M^{[3]}\frac{r_M^{[3]}}{r_{pm}^3}\sin\left(\theta^{[3]}-L_M^{[3]}\right), \vspace{5pt}\\

a_{pr}^{[4]} & = & \displaystyle -\mu_S^{[4]}\left[\frac{r^{[4]}}{r_{psM}^3} + \left(\frac{r_M^{[4]}}{r_{psM}^3}-\frac{1}{\left(r_M^{[4]}\right)^2}\right)\cos\left(\theta^{[4]}-L_M^{[4]}\right)\right], \vspace{5pt}\\

a_{p\theta}^{[4]} & = & \displaystyle - \mu_S^{[4]}\left(-\frac{r_M^{[4]}}{r_{psM}^3}+\frac{1}{\left(r_M^{[4]}\right)^2}\right)\sin\left(\theta^{[4]}-L_M^{[4]}\right),
\end{array} 
\end{equation}
where
\begin{equation}\label{DummyVariablesForDistance}
\begin{array}{lcl}
r_{psE} & = & \displaystyle \sqrt{\left(r_E^{[2]}\right)^2 + \left(r^{[2]}\right)^2 + 2r^{[2]}r_E^{[2]}\cos\left(\theta^{[2]} - L_E^{[2]}\right)}, \vspace{5pt}\\

r_{pe} & = & \displaystyle \sqrt{\left(r_E^{[3]}\right)^2 + \left(r^{[3]}\right)^2 - 2r^{[3]}r_E^{[3]}\cos\left(\theta^{[3]} - L_E^{[3]}\right)}, \vspace{5pt}\\

r_{pm} & = & \displaystyle \sqrt{\left(r_M^{[3]}\right)^2 + \left(r^{[3]}\right)^2 - 2r^{[3]}r_M^{[3]}\cos\left(\theta^{[3]} - L_M^{[3]}\right)}, \vspace{5pt}\\

r_{psM} & = & \displaystyle \sqrt{\left(r_M^{[4]}\right)^2 + \left(r^{[4]}\right)^2 + 2r^{[4]}r_M^{[4]}\cos\left(\theta^{[4]} - L_M^{[4]}\right)}. \vspace{5pt}\\

\end{array}
\end{equation}

\end{document}